\documentclass[12pt]{amsart} 
\usepackage{amssymb}
\usepackage{latexsym}
\usepackage{amsmath}
\usepackage{amsthm}

\usepackage{doi}

\DeclareFontFamily{U}{mathx}{}
\DeclareFontShape{U}{mathx}{m}{n}{<-> mathx10}{}
\DeclareSymbolFont{mathx}{U}{mathx}{m}{n}
\DeclareMathAccent{\widecheck}{0}{mathx}{"71}

\usepackage{tikz-cd} 

\newtheorem{theorem}{Theorem}[section]

\newtheorem{lemma}[theorem]{Lemma}
\newtheorem{corollary}[theorem]{Corollary}
\newtheorem{proposition}[theorem]{Proposition}

\theoremstyle{definition}

\theoremstyle{remark}

\numberwithin{equation}{section}
\numberwithin{theorem}{section}
\allowdisplaybreaks
\usepackage{color}

\newcommand{\R}{\mathbb R}

\newcommand{\C}{\mathbb C}

\newcommand{\logl}{L\textnormal{log} L}
\newcommand{\logld}{L(\textnormal{log} L)^2}
\newcommand{\lexp}{L_{\textnormal{exp}}}
\newcommand{\logla}{L(\textnormal{log } L)^{\alpha}}
\newcommand{\loglb}{L(\textnormal{log } L)^{\beta}}
\newcommand{\tlog}{T_{\textnormal{log}}}

\title[The finite Hilbert transform]
{The finite Hilbert transform on $(-1,1)$}

\author[G. P.  Curbera]{Guillermo P. Curbera}
\address{Facultad de Matem\'aticas \& IMUS,
Universidad de Sevilla, 
Calle Tarfia s/n,  Sevilla 41012, Spain}
\email{curbera@us.es}

\author[S. Okada]{Susumu Okada}
\address{112 Marcorni Crescent, Kambah, ACT 2902, Australia}
\email{susbobby@grapevine.com.au}

\author[W.J. Ricker]{Werner J. Ricker}
\address{Math.--Geogr.\  Fakult\"at, Katholische Universit\"at
Eichst\"att--Ingolstadt, D--85072 Eichst\"att, Germany}
\email{werner.ricker@ku.de}

\thanks{The first author acknowledges the support of PID2021-124332NB-C21
(FEDER(EU)/Ministerio de Ciencia e Innovaci\'on-Agencia Estatal de Investigaci\'on)
and FQM-262 (Junta de Andaluc\'{\i}a).}

\date{\today}


\begin{document}


\begin{abstract} 
We present a detailed survey of  recent developments in the study
of the finite Hilbert transform  and its corresponding 
inversion problem in rearrangement invariant  spaces on $(-1,1)$.
\end{abstract}

\maketitle

\tableofcontents


\section{Introduction}
\label{S1}


Given $f\in L^1(-1,1)$ its \textit{finite Hilbert transform} $T(f)$   is the principal value integral
\begin{equation}\label{T}
T(f)(t)=\lim_{\varepsilon\to0^+} \frac{1}{\pi}
\left(\int_{-1}^{t-\varepsilon}+\int_{t+\varepsilon}^1\right) \frac{f(x)}{x-t}\,dx, 
\end{equation}
which exists for  a.e.\   $t\in(-1,1)$ and is a measurable function.
Throughout the paper we will shorten the terminology to FHT. Its study is
intimately related to the solution of the \textit{airfoil equation}, that is, 
for $g$ a suitable given function, find all functions $f$ which satisfy
\begin{equation}\label{airfoil}
g(t)=\mathrm{p.v.} \frac{1}{\pi} \int_{-1}^{1}\frac{f(x)}{x-t}\,dx,
\quad \mathrm{a.e. }\; t\in(-1,1).
\end{equation}

In the early times  of Aerodynamics the study and resolution of the 
airfoil  equation played a central role: 
\begin{quotation}
``In the theory of the two dimensional flow of an ideal fluid past a thin airfoil there arise 
two types of problem, which may be called the ``thickness'' problem and 
the ``lifting'' problem; they lead to two different types of boundary value 
problem for the complex velocity $w=u-iv$ in the complex plane of flow, $z=x+iy$,''
\end{quotation}
\cite[p.357]{ChR}. The study of the second of these two problems 
lead to the \textit{``integral equation of the 
lifting problem''}, namely
\begin{equation*}
v_0(x)=\frac{1}{\pi}
\int_{\alpha_1}^{\alpha_2} \frac{u_0(s)}{s-x}\,ds,
\end{equation*}
where $v_0(x):=v(x,+0)$ and $u_0(x):=u(x,+0)$;
\cite[p.358-359]{ChR}. The early treatment of the lifting problem and its inversion 
began in the 1920s with the works of Betz, Birnbaum, Carleman, Glauert, and  Munk,  
and continued  in the 1930s with the work of Hamel and S\"ohngen; for a detailed account
see  \cite[\S1]{ChR} and \cite{So-1}. 
The study of the $L^p$-theory for the FHT began in the 1950s 
with the work of Tricomi,  \cite{T1, T}, and S\"ohngen, \cite{So-2}, and later continued 
by Widom, \cite{W},  J\"orgens, \cite{J}, and others. 
The 1991 paper of Okada and Elliott, \cite{OE}, gave a
complete and compact presentation of these results (with alternative
proofs) together with  a clear account of the state-of-the-art of the $L^p$-theory at that time.

Simultaneously,  the consideration of 
problems arising in mathematical physics, in particular,  in elasticity theory,
led the Soviet School to  study general one-dimensional singular integral operators 
closely related to the FHT and defined on curves more general than $(-1,1)$. 
In this direction  the works of
Duduchava,  \cite{D}, Gakhov,  \cite{G}, Gohberg and Krupnik, \cite{GK-1,GK-2}, 
Khvedelidze, \cite{Kh}, 
Mikhlin and Pr\"ossdorf, \cite{MP}, and Muskhelishvili, \cite{Mu}, amongst others,
should be  highlighted; see also the references therein.
This important topic is beyond the scope of the present article, in which 
the FHT is treated exclusively on the particular curve $(-1,1)$. 
We will follow what King called ``the Tricomi approach'', \cite[\S11.4]{K}.


Our central aim is to discuss a series of six recent papers, 
\cite{COR-ampa,COR-qm,COR-mn,COR-am,COR-asnsp,COR-prep}, 
where the focus of studying the FHT was transferred from the traditional family
of spaces $L^p(-1,1)$, for $1<p<\infty$, to the significantly larger class
of rearrangement invariant (in brief, r.i.) spaces $X$ over $(-1,1)$. 
That this class of spaces  is the most  suitable one  to consider 
is illustrated by the fact that  $T\colon X\to X$ is injective if 
and only if $L^{2,\infty}(-1,1)\not\subseteq X$ 
and (for $X$ separable) that $T\colon X\to X$ has a 
non-dense range if and only if $X \subseteq L^{2,1}(-1,1)$. 
Here, the Lorentz spaces $L^{2,1}(-1,1)$ and $L^{2,\infty}(-1,1)$ are 
r.i.\ spaces.

It is appropriate that we begin by recalling the $L^p$-theory of the FHT, 
with the intention of providing the basis for a better understanding
of the recent results. 
This is done in Section \ref{S3}.  The extension of the results from the 
$L^p$-setting   to
r.i.\ spaces is presented in Sections \ref{S4} and \ref{S5}. In the event that $X$ has non-trivial
Boyd indices, there is a close connection between the $L^p$-theory and  theory of r.i.-spaces
(see Lemma \ref{l-4.1}).
In Section \ref{S6} it is shown that the natural domain of the FHT operator
$T\colon X\to X$ ``is actually $X$ itself'', meaning that $T$ cannot be extended to any larger 
domain space (while still maintaining its values in $X$).
The FHT does not map the classical Zygmund space $\logl$ into itself. However, it does 
map $\logl$ continuously into $L^1$. This operator is investigated in Section \ref{S7},
where extrapolation plays an important role.
Widom determined completely 
the spectrum of $T\colon L^p(-1,1)\to L^p(-1,1)$, for
all $1<p<\infty$. Section \ref{S8} is devoted to extending these results to
$T\colon X\to X$ for  r.i.\ spaces $X$ with non-trivial Boyd indices.
Section \ref{S9} discusses the use of the theory of integration with respect
to Banach space-valued measures to provide an integral representation of 
$T\colon X\to X$.

It should be pointed out that  results on the FHT have recently found applications 
to problems arising in image reconstruction; see, for example, \cite{KT}, \cite{Si}, and the 
references therein.

Note that the definition of the FHT presented in  \eqref{T}  
coincides with the one used by  Tricomi, \cite[\S4.3]{T}.
Other  definitions of the FHT also appear in the literature, differing
from the one above by a multiplicative constant: King uses $-T$, \cite[Ch.11]{K}, 
whereas Widom, \cite{W}, and J\"orgens, \cite[\S13.6]{J}, use $T/i$.
In Section \ref{S8} we will use $T/i$ in order to be consistent with the presentation of the results
by Widom and J\"orgens.


\section{Preliminaries}
\label{S2}

The setting of this paper is the measure space consisting of $(-1,1)$ 
equipped with its Borel $\sigma$-algebra $\mathcal{B}$ and  Lebesgue measure  $m$
on $\R$ restricted to $\mathcal{B}$.
We  denote by $\text{sim }\mathcal{B}$ the vector space of 
all $\mathbb{C}$-valued, $\mathcal{B}$-simple functions on $(-1,1)$  and by
$L^0(-1,1)=L^0$ the space (of equivalence classes) of all $\mathbb{C}$-valued
measurable functions, endowed with the topology of convergence in measure.
The space $L^p(-1,1)$ is denoted simply by $L^p$, for $1\le p\le\infty$.

A \textit{Banach function space} (B.f.s.) $X$ on  $(-1,1)$ is a
Banach space  $X\subseteq L^0$ satisfying
the ideal property, that is, $g\in X$ and $\|g\|_X\le\|f\|_X$
whenever $f\in X$ and $|g|\le|f|$ a.e.  
The \textit{associate space} $X'$  of $X$ is the B.f.s.\  which consists  of all
functions $g\in L^0$ satisfying $\int_{-1}^1|fg|<\infty$, for every
$f\in X$, equipped with the norm
$\|g\|_{X'}:=\sup\{|\int_{-1}^1fg|:\|f\|_X\le1\}$. 
The space $X'$ is isometrically isomorphic to 
a closed subspace of the dual Banach space $X^*$ of $X$. 
Moreover, if $f\in X$ and $g\in X'$, then $fg\in L^1$ and
$\|fg\|_{L^1}\le \|f\|_X \|g\|_{X'}$, i.e., H\"older's inequality
is available. The second  associate space $X''$ 
of $X$ is defined to be $(X')'$. The norm in $X$ is said to be
absolutely continuous (in brief, a.c.) if,
for every $f\in X$, we have $\|f\chi_A\|_X\to0$ whenever $m(A)\to0$.
In case the norm in $X$ is not a.c., we can consider the 
closed subspace $X_a$ consisting of all the absolutely continuous elements of $X$, that is,
of all $f\in X$ such that $\|f\chi_A\|_X\to0$ whenever $m(A)\to0$.
The space $X$ satisfies the Fatou property  if, whenever  $\{f_n\}_{n=1}^\infty\subseteq X$ satisfies
$0\le f_n\le f_{n+1}\uparrow f$ a.e.\ with $\sup_n\|f_n\|_X<\infty$,
then $f\in X$ and $\|f_n\|_X\to\|f\|_X$.   
As in  \cite{BS},  
\textit{all} B.f.s.' $X$ (hence, all r.i.\ spaces)   
are assumed to satisfy the Fatou property. In this case $X''=X$ and hence,
$f\in X$ if and only if $\int_{-1}^1|fg|<\infty$, for every $g\in X'$.
Moreover, $X'$ is a norm-fundamental subspace of $X^*$, that is, 
$\|f\|_X=\sup_{\|g\|_{X'}\le1} |\int_{-1}^1fg|$ for $f\in X$, \cite[pp.12-13]{BS}.
If $X$ is separable, then $X'=X^*$.

A \textit{rearrangement invariant} (r.i.) space $X$ on $(-1,1)$ is a
B.f.s.\  having the property that whenever $g^*\le f^*$ with $f\in X$,  
then $g\in X$ and $\|g\|_X\le\|f\|_X$.
Here $f^*\colon[0,2]\to[0,\infty]$ is 
the decreasing rearrangement of $f$, that is, the
right continuous inverse of its distribution function:
$\lambda\mapsto m(\{t\in (-1,1):\,|f(t)|>\lambda\})$ for $\lambda\ge0$.
The associate space $X'$ of a r.i.\ space $X$ is again a r.i.\ space.
Every r.i.\ space $X$  satisfies 
$L^\infty\subseteq X\subseteq L^1$ with continuous inclusions. 
The fundamental function of $X$ is defined by 
$\varphi_X(t):=\|\chi_{A}\|_X$ for $A\in\mathcal{B}$ 
with $m(A)=t$, for $t\in[0,2]$.

The family of r.i.\ spaces includes many classical spaces 
appearing in analysis, such as  the Lorentz $L^{p,q}$ spaces, 
\cite[Definition IV.4.1]{BS}, Orlicz $L^\varphi$ spaces 
\cite[\S4.8]{BS}, Marcinkiewicz $M_\varphi$ spaces, 
\cite[Definition II.5.7]{BS}, Lorentz $\Lambda_\varphi$ spaces, \cite[Definition II.5.12]{BS},
and the Zygmund $L^p(\text{log L})^\alpha$ spaces, 
\cite[Definition IV.6.11]{BS}. In particular,  $L^p=L^{p,p}$, for $1\le p\le \infty$.  
The space weak-$L^1$, denoted by $L^{1,\infty}(-1,1)=L^{1,\infty}$,  
is a quasi Banach space, \cite[Definition IV.4.1]{BS}, and satisfies
$L^1\subseteq L^{1,\infty} \subseteq L^0$, with both inclusions continuous.

An important role is played by the Marcinkiwiecz space 
$L^{2,\infty}(-1,1)=L^{2,\infty}$, also known as  weak-$L^2$, \cite[Definition IV.4.1]{BS}. 
It consists of those functions $f\in L^0$  satisfying
$$
f^*(t)\le \frac{M_f}{t^{1/2}},\quad 0<t\le2,
$$
for some constant $M_f>0$. Consider the function 
\begin{equation}\label{w}
 w(x) : = \sqrt{1-x^2}, \quad x \in (-1,1),
\end{equation}
which pervades the theory of the FHT.
Since the  decreasing rearrangement  of the function $1/w$ on  $(-1,1)$
is the function $t\mapsto 2/t^{1/2}$ on $(0,2)$,  it follows that $1/w$ belongs to
$L^{2,\infty}$. Actually, for any r.i.\ space $X$  it is the case that 
$1/w\in X$ if and only if $L^{2,\infty}\subseteq X$. 
Consequently, $L^{2,\infty}$ is the \textit{smallest} r.i.\ space  which contains  $1/w$.

Standard references concerning B.f.s.' and r.i.\ spaces 
are \cite{BS}, \cite{KPS}, \cite{LT}.


Let $Y$ be a  Banach space with norm $\|\cdot\|_Y$  and   dual space $Y^*$,
equipped with the usual  dual norm $\|\cdot\|_{Y^*}$. The identity operator on 
$Y$ is denoted by $I_Y$.  Let $B(Y)$  denote  the vector space of all 
continuous linear operators from $Y$ into itself.    Given $S \in B(Y)$, 
denote by $S^* \in B(Y^*)$  its corresponding adjoint operator.   
By $\mathrm{Ker}(S)$ and   $\text{R} (S)$  we denote the kernel 
and the range space of $S$, respectively, that is, $\mathrm{Ker}(S):= S^{-1}(\{0\})$ 
and $ \text{R} (S) :=\{S(y): y \in Y\}$.
It is known  that $\mathrm{Ker}(S^*) $ equals the annihilator 
$ \text{R}(S)^\perp$ of $R(S)$, defined to be the closed linear subspace of $Y^*$ 
consisting of all functionals $y^* \in Y^*$ which vanish on $\text{R}(S)$.  
The dual space $\big( Y/\overline{\text{R}(S)}\big) ^*$ of the quotient Banach space 
$ Y/\overline{\text{R}(S)}$ is  isometrically isometric to  $\text{R}(S)^\perp$, where
$\overline{\text{R}(S)}$ denotes the closure of $\text{R}(S)$  in $Y$.

An operator $S \in B(Y)$ with closed range   is  called a \textit{Fredholm operator}  
if   $\dim(\mathrm{Ker}(S))< \infty$ and 
$\dim(Y/ \text{R}(S))< \infty$. In this case there exist projections 
$P,\,Q \in B(Y)$  satisfying  both $\text{R}(P) = \mathrm{Ker}(S)$ and 
$\mathrm{Ker}(Q) =\text{R}(S)$. Furthermore, 
there exists a unique operator $R \in B(Y)$ satisfying
$$
RS= I_Y-P, \quad SR = I_Y-Q,\quad PR = 0\quad \text{and} \quad RQ = 0;
$$
see \cite[Theorem 5.4]{J}, for example.
The operator $R$   is said to be the \textit{pseudo-inverse} of $S$ relative to 
the projections $P$ and $Q$ and the integer 
$ \kappa(S): = \dim (\mathrm{Ker}(S)) - \dim (Y/ \text{R}(S))$
is called the \textit{index} of the Fredholm operator $S$.


\section{The airfoil equation: $L^p$-theory}
\label{S3}

In this section we sketch some relevant aspects of the 
$L^p$-theory for the FHT and the airfoil equation.
A celebrated theorem of M. Riesz  
states, for each $1<p<\infty$, that
the Hilbert transform operator  $H$ maps $L^p(\R)$ continuously  into 
itself, \cite[Theorem III.4.9(a)]{BS}. Since the FHT  given
in \eqref{T} can be written as  $Tf=\chi_{(-1,1)}H(f\chi_{(-1,1)})$, 
it follows  that the linear operator $f\mapsto T(f)$
maps $L^p$ continuously  into itself. We will denote this operator by $T_p$.
However,   $T$ is not continuous on $L^\infty$ nor   on $L^1$.

The following two formulae  are fundamental  for the study of the FHT and its inversion.
\begin{itemize}
\item[(a)] 
The \textit{Parserval formula} holds for a pair of functions $f, g\in L^1$ means that 
both integrands belong to $L^1$ and
\begin{equation}\label{parseval}
\int_{-1}^{1}fT(g)=-\int_{-1}^{1}gT(f).
\end{equation}
\item[(b)] The \textit{Poincar\'e-Bertrand formula} holds  for a pair 
of functions $f, g\in L^1$, with all the terms
finite a.e., means that
\begin{equation}\label{poincare}
T(gT(f)+fT(g))=T(f)T(g)-fg,\quad \text{a.e. on }  (-1,1).
\end{equation}
\end{itemize}


In order to solve the airfoil equation \eqref{airfoil} and find the \textit{inversion formula}
for the FHT, Tricomi   argued in \cite[\S4]{T1} in the following way.
Let $g$ be given. Denote by $f$ the solution (if it exists) of the equation $T(f)=g$.
Applying the Poincar\'e-Bertrand formula \eqref{poincare} to the pair of functions
$f$ and   $w$ (cf.\   \eqref{w}), and noting   that $T(w)(t)=-t$, \cite[(11.57)]{K}, yields
\begin{equation}\label{3.3}
T\big(-xf(x)+w(x)g(x)\big)(t)=-tg(t)-w(t)f(t).
\end{equation}
Direct computation shows that
\begin{equation*}
T\big(xf(x)\big)(t)=T\big((x-t+t)f(x)\big)(t)=\frac{1}{\pi} \int_{-1}^{1}f(x)\,dx+ tT(f(x))(t),
\end{equation*}
which, in view of \eqref{3.3}, implies that
\begin{equation*}
w(t)f(t)=t\big(T(f(x))(t)-g(t)\big)-T\big(w(x)g(x)\big)(t)+\frac{1}{\pi} \int_{-1}^{1}f(x)\,dx.
\end{equation*}
Setting $C:=\frac{1}{\pi} \int_{-1}^{1}f(x)\,dx$, we arrive at
$$
f=\frac{-1}{w} T(wg)+\frac{C}{w}.
$$
These, and related computations, reveal the need for a deep 
analysis of the integrability properties of the functions involved.

Regarding   Parseval's formula, given $1 < p < \infty$ and 
its conjugate index $p'$, i.e., $1/p + 1/{p'} = 1$, the identity \eqref{parseval}
is valid for every pair $ f \in L^p$ and $g \in L^{p'}$,  \cite[\S4.3 (2)]{T}.
An immediate consequence   is that the adjoint operator $(T_p)^*$ of $T_p$ 
is given by $(T_p)^* = -T_{p'}$.

Regarding the Poincar\'e-Bertrand formula \eqref{poincare}, 
it was  proved by Tricomi  for a pair of functions  $f\in L^p$ and $g\in L^q$ 
whenever the indices satisfy $1/p+ 1/q <1$, \cite[\S4.3 (4)]{T}).  
In 1977, Love  established that \eqref{poincare} also holds for 
all $ f \in L^p$ and $g \in L^{p'}$, \cite[Corollary]{L}.  
A proof can also be obtained by using   Chebyshev polynomials, 
\cite[Theorem 2.7]{O1}. For the earlier history of the Poincar\'e-Bertrand formula, 
we refer to  \cite[\S 2.13 \& \S 4.23]{K} and to the Introduction in \cite{L}.


Considerations concerning weighted versions of the FHT are important.
In this regard the following result of Khvedelidze, \cite{Kh}, is relevant. Further proofs 
of it occur in \cite[Lemma I.4.2]{GK-1} and in \cite[Theorem II.3.1]{MP}.

\begin{theorem} \label{t-3.1} 
Let  $1 < p < \infty$ and $\rho$ be the weight function
$$
\rho(x) : = (1-x)^\gamma (1+x)^\delta, \quad x \in (-1,1),
$$
where $\gamma, \delta \in (-1/p, 1/p')$. Then the function 
$\rho T(f/\rho)$ belongs to $L^p$ for every $f \in L^p$ 
and the resulting linear operator
$$
f \mapsto \rho T(f/\rho), \quad f \in L^p,
$$
is continuous from $L^p$ into $L^p$.
\end{theorem}

The particular weight function $ w$ defined in \eqref{w} 
plays a fundamental role in the study of the FHT since $T(1/w) = 0$, \cite[\S4.3 (7)]{T},
and because it determines the kernel of $T$, \cite[\S4.3 (14)]{T}, namely,
\begin{equation}\label{3.4}
\text{$f \in \bigcup_{1< p < \infty} L^p$ satisfies $T(f)= 0$ $\iff$  
$f = C/w$ for some $C \in \C$.}
\end{equation}
Since $1/w$ belongs to $L^p  \setminus L^2$, for every $1< p< 2$, 
it follows that
 \begin{equation}\label{3.5}
\mathrm{Ker}(T_p) = \mathrm{span}\{1/w\},\;\;  1 < p < 2,
\qquad  \mathrm{Ker}(T_p) = \{0\},\;\; 2\le p < \infty.
\end{equation}
%


For each $f\in L^1$, the function $\widehat T(f)\in L^0$ is 
defined pointwise a.e.\  in $(-1,1)$ by
\begin{equation}\label{3.6}
\widehat{T}(f):=- \frac{1}{w}T(fw) .
\end{equation}


Fix $1<p< 2$. Then the conditions of 
Theorem \ref{t-3.1} are satisfied for  $\gamma=\delta=-1/2$ (with $\rho: = 1/w$).   
It follows that the restriction of $\widehat T$ from $L^1$ to $L^p$
defines a continuous linear operator $\widehat{T}_p\colon L^p\to L^p$, namely,
$$
\widehat{T_p}(f) :=- \frac{1}{w}T_p(fw), \quad f \in L^p.
$$
Moreover, as $1/w\in L^p$, the linear operator  $P_p\colon L^p \to L^p$ defined by
$$
P_p(f) : = \left(\frac{1}{\pi}\int_{-1}^1 f(x)\,dx\right)\frac{1}{w}, \quad f \in L^p,
$$
is a continuous  projection onto the one-dimensional linear 
subspace $\mathrm{span}\{1/w\}$ of $L^p$.  In particular, 
\eqref{3.5} shows that $\mathrm{Ker}(T_p)= \text{R}(P_p)$.  
The operator $\widehat{T}_p$ turns out to be  the pseudo-inverse 
of $T_p$ relative to the projection $P_p$ and  the zero operator, 
as  formulated  in the following theorem, \cite[Proposition 2.4]{OE};  see also its proof.


\begin{theorem}\label{t-3.2}
For $T_p\colon L^p\to L^p$ with $1<p<2$ the following statements hold.
\begin{itemize}
\item[(i)] $\mathrm{Ker}(T_p)=\mathrm{span}\{1/w\}$.
\item[(ii)] The continuous linear operator $\widehat{T}_p\colon L^p\to L^p$ 
satisfies $ T_p \widehat{T}_p=I_{L^p}$ and 
$$
\int_{-1}^1 \widehat{T}_p (f)(x)\,dx = 0,\quad f \in L^p.
$$
\item[(iii)] The operator $T_p\colon L^p\to L^p$ is surjective.
\item[(iv)] The identity $\widehat{T}_p T_p=I_{L^p}-P_p$ holds.
\item[(v)] The operator $\widehat{T}_p$ is  an isomorphism onto its range 
$\text{R}(\widehat{T}_p)$, which is given by
$$
\text{R}(\widehat{T}_p)=\left\{f\in L^p: \int_{-1}^{1} f(x)\,dx=0\right\}.
$$
\item[(vi)] The following direct sum decomposition of $L^p$ holds:
$$
L^p =\text{R}(\widehat{T}_p) \oplus \mathrm{span}\{1/w\}.
$$
\item[(vii)] The operator $T_p$ is  Fredholm  with index $\kappa(T_p) =1$ 
and $\widehat{T}_p$ is its pseudo-inverse  
relative to the projection $P_p$ and the zero operator.
\end{itemize}
\end{theorem}


For $f\in L^1$ satisfying $f/w\in L^1$, the function $\widecheck T(f)\in L^0$ 
is defined pointwise a.e.\ in $(-1,1)$ by
\begin{equation}\label{3.7}
\widecheck{T}(f) :=-w\,T\,\Big(\frac{f}{w}\Big) .
\end{equation}


For $2 < p < \infty$   we are in the setting of 
Theorem \ref{t-3.1}  for $\gamma=\delta=1/2$  and $\rho: = w$. 
It follows that the restriction of $\widecheck T$  to $L^p$
defines a continuous linear operator 
$\widecheck{T}_p\colon L^p\to L^p$, namely,
$$
\widecheck{T}_p(f) :=-w\,T_p\,\Big(\frac{f}{w}\Big), \quad f \in L^p.
$$
Furthermore, the linear operator  $Q_p\colon L^p \to L^p$ defined by
$$
Q_p(f):=\left(\frac1\pi\int_{-1}^1
\frac{f(x)}{w(x)}\,dx\right) \mathbf{1},\quad f\in L^p,
$$
where $\mathbf{1}:=\chi_{(-1,1)}$, is a continuous  
projection onto the one-dimensional linear subspace   $\mathrm{span}\{\mathbf{1}\}$.  
In particular, $\text{R}(T_p)=\mathrm{Ker}(Q_{p})$.  
The operator $\widecheck{T}_p$ turns out to be  the pseudo-inverse 
of $T_p$ relative to the projection $Q_p$ and  the zero operator, 
as  stated in the following result, \cite[Proposition 2.6]{OE}; see also its proof.


\begin{theorem}\label{t-3.3}
For $T_p\colon L^p\to L^p$ with $2<p<\infty$ the following statements hold.
\begin{itemize}
\item[(i)] The operator $T_p\colon L^p\to L^p$ is injective.
\item[(ii)] The continuous  linear operator $\widecheck{T}_p\colon L^p\to L^p$ 
satisfies $\widecheck{T_p} T_p=I_{L^p}$.
\item[(iii)] The identity $T_p\widecheck T_p=I_{L^p} -Q_p$ holds in $L^p$.
\item[(iv)] The range of $T_p$ is the closed subspace of $L^p$ given by
$$
\text{R}(T_p)=\left\{f\in L^p: \int_{-1}^{1} \frac{f(x)}{w(x)}\,dx=0\right\}
=\mathrm{Ker}(Q_{p}).
$$
Moreover,  $\widecheck T_p$ is an isomorphism from $\text{R}(T_p)$ onto $L^p$.
\item[(v)] The following direct sum decomposition of $L^p$ holds:
\begin{equation}\nonumber
L^p = \text{R}(T_p)\oplus\mathrm{span}\{\mathbf{1}\}.
\end{equation}
\item[(vi)]  The operator $T_p$ is  Fredholm  
with $\kappa(T_p) =-1$ and $\widecheck{T}_p$ is its pseudo-inverse  
relative to the  zero operator  and the  projection $Q_p$.
\end{itemize}
\end{theorem}


Theorems \ref{t-3.2} and \ref{t-3.3} lead directly  to  the 
inversion formula  for solving the airfoil equation \eqref{airfoil}
within $L^p$  whenever $1< p < 2$ and $2< p< \infty$;
see Corollaries 2.5 and 2.8 in \cite{OE}.

\begin{corollary}\label{c-3.4}  
The following inversion formulae hold.
\begin{itemize}
\item[(i)]  Let  $1< p < 2$. Given $g \in L^p$, a function 
$ f \in L^p$ is a solution of the airfoil equation \eqref{airfoil}   if and only if
$$
f = -\frac{1}{w}\,T(g w) + \frac{C}{w},
$$
for a constant $C\in\C$, in which case $C=(1/\pi)\int_{-1}^1 f(x)\,dx $. 
\item[(ii)]   Let $2< p < \infty$.  Let  $g\in L^p$ satisfy  
$\int_{-1}^{1} g(x)/w(x)\,dx = 0 $.  Then,  the airfoil equation \eqref{airfoil} 
admits a unique solution $f\in L^p$ given by
$$
f= -wT\Big(\frac{g}{w}\Big).
$$
\end{itemize}
\end{corollary}


We end this section with some historical comments. 
The arguments used to establish Theorems \ref{t-3.2} and \ref{t-3.3} and Corollary \ref{c-3.4}
above are built on Tricomi's work in \cite[\S 4.3]{T}.  
The main tools are the Parseval identity \eqref{parseval}, the Poincar\'e-Bertrand 
identity \eqref{poincare} and Khvedelidze's Theorem  \ref{t-3.1}.  
It should be noted that   \eqref{poincare} was only available 
to S\"ohngen and Tricomi under more restrictive conditions on a pair 
of functions $f, g$   and that Theorem \ref{t-3.1} was unknown  in the early 1950s.  
A complete and compact presentation of these results (with alternative
proofs) was given in \cite{OE}. 
The inversion formula in Corollary \ref{c-3.4}(i),  with certain restrictions, appeared  in  
\cite[(22)]{So-1}  and again in \cite[\S3,4]{T1}  and  \cite[Satz 6]{So-2}   
with  less stringent  conditions.  The same inversion formula is also presented  in 
\cite[\S4.3 (16)]{T}  together with a further explanation in the footnote on p.179. 
An alternative proof of Theorem  \ref{t-3.3} can be
obtained  by applying Theorem \ref{t-3.2} with $p'$ in place of $p$.

For the case  $p=2$,  the operator $T_2\colon L^2\to L^2$ is injective (by 
\eqref{3.5}) and  its adjoint equals  $-T_2$ (by   Parseval's formula) and 
hence, it  is also injective. This implies that the range  of $T_2$ is 
a proper, dense, linear subspace of $L^2$. Note that the function
$f(x):= x/w(x)$, for $x \in (-1,1)$, satisfies
$f\in L^p\setminus L^2$ for each $1< p < 2$ but,  $T(f) =\mathbf{1} \in L^2$
(for further interesting examples, see \cite[Lemma 4.3 \& Note 4.4]{OE}).   
So, $T_2$  is not a Fredholm operator and its inverse $T_2^{-1}$ is 
an unbounded operator with dense domain $\mathcal{R}(T_2)$.  
For a detailed study and further  properties of $T_2$ and its inversion formulae, see
\cite[Sections 3 \& 4]{OE}.


The study of singular integral operators in  $L^2$ 
which are not Fredholm   has the difficulty 
that their inverse is an unbounded operator.  
This problem was  already noted  for the FHT in \cite[p.44]{So-2}. 
To deal with this feature,  the FHT has also been considered
as acting on certain weighted $L^2$-spaces. This has been
undertaken for the weights $w$, $1/w$, $\sigma$ and $1/\sigma$,
where $\sigma(x) :=(1-x)^{-1/2} (1+x)^{1/2}$  for $x \in (-1,1)$.
For each of these weights the   FHT is a Fredholm operator on the 
corresponding weighted $L^2$-space, a fact which has been 
demonstrated to be useful in solving certain singular integral 
equations, \cite{Sch1,Sch2}. Motivated by \cite{Sch1,  Sch2}, 
the paper \cite[\S 4]{O2} studies the FHT on general 
weighted $L^p$-spaces with $1< p < \infty$. 
For further investigations of the FHT on   weighted  $L^p$-spaces 
and related  Sobolev-type spaces 
we refer to  \cite{APS, BHS1, BHS2, OP} (and the references therein).  
In \cite{EO}  a further class of  weighted 
Sobolev spaces was introduced in which  the FHT 
turns out to be  a continuous operator.


\section{The finite Hilbert transform  in rearrangement invariant spaces}
\label{S4}

A celebrated result of Boyd leads to the extension of  the classical result of M. Riesz,
asserting the continuity of the Hilbert transform $H$ on  
$L^p(\R)$, for  $1<p<\infty$, to a larger class of r.i.\ spaces. 
Given any r.i.\ space $Y$ over $\R$, Boyd associated two indices, 
$\underline{\alpha}_Y$ and  $\overline{\alpha}_Y$, to $Y$ which satisfy 
$0\le\underline{\alpha}_Y\le \overline{\alpha}_Y\le1$,
and proved that  $H$ acts continuously on
$Y$  if and only if those  indices  are non-trivial, that is, 
$0<\underline{\alpha}_Y\le \overline{\alpha}_Y<1$,  \cite[Theorem III.5.18]{BS}.
The indices $\underline{\alpha}_Y$ and  $\overline{\alpha}_Y$ are called 
the lower and upper Boyd indices  of $Y$, respectively.
Regarding the FHT acting on  r.i.\ spaces on $(-1,1)$, 
the analogous characterization as above is valid:  $T$ acts continuously on
a r.i.\ space $X$ on $(-1,1)$ if and only if $X$ has non-trivial Boyd indices; 
see, for example, \cite[pp.170-171]{KPS}.


The construction of the Boyd indices in the case of a r.i.\ space $X$ on $(-1,1)$ 
proceeds as follows
(for the  setting of more general measure spaces see \cite[\S III.5]{BS}).  
Given such a r.i.\ space  $X$, the Luxemburg representation theorem 
ensures that there exists another  r.i.\  space  $\widetilde X$ on $(0,2)$
such that $\|f\|_X=\|f^*\|_{\widetilde X}$ for $f\in X$, \cite[Theorem II.4.10]{BS}.
The dilation operator $E_t$ for $t>0$ is defined, for 
each $f\in \widetilde X$, by $E_t(f)(s):=f(st)$ for $0\le s\le \min\{2,1/t\}$ and $E_t(f)(s)=0$ 
for  $\min\{2,1/t\}< s\le 2$. The operator $E_t\colon \widetilde X\to \widetilde X$  is continuous 
with $\|E_{1/t}\|_{\widetilde X\to \widetilde X}\le \max\{t,1\}$. The \textit{lower} and \textit{upper 
Boyd indices} of  $X$ are then defined, respectively, by
$$
\underline{\alpha}_X\,:=\,\sup_{0<t<1}\frac{\log \|E_{1/t}\|_{\widetilde X\to \widetilde X}}{\log t}
\;\;\mbox{and}\;\;
\overline{\alpha}_X\,:=\,\inf_{t>1}\frac{\log \|E_{1/t}\|_{\widetilde X\to \widetilde X}}{\log t} ;
$$
see  \cite[Definition III.5.12]{BS}. There are other indices that will be needed 
when studying the spectrum of $T$ in Section \ref{S8}. 
The lower and upper fundamental indices, 
$\underline{\beta}_X$ and $\overline{\beta}_X$,  are  defined by
$$
\underline{\beta}_X:= \sup_{0<t<1}\frac{\log M_{\varphi_X}(t)}{\log t}
\;\;\mbox{and}\;\;
\overline{\beta}_X:= \inf_{t>1}\frac{\log M_{\varphi_X}(t)}{\log t},
$$
\cite[pp. 177-178]{BS}, respectively, where  
$$
M_{\varphi_X}(t):=\sup_{0< s < \min\{2,2/t\}}\frac{\varphi_X(st)}{\varphi_X(s)},
\quad t\in(0,\infty),
$$
and $\varphi_X$ is the fundamental function of $X$ (see Section \ref{S2}). 
The following inequalities hold:
$$
0\le\underline{\alpha}_X\le \underline{\beta}_X\le
\overline{\beta}_X\le \overline{\alpha}_X\le 1.
$$
For $X=L^p$ with $1<p<\infty$, it is known that 
$\underline{\alpha}_{L^p}= \overline{\alpha}_{L^p}=1/p$.

The class of r.i.\ spaces on $(-1,1)$ with non-trivial Boyd indices
is closely connected to the family of $L^p$-spaces for $1<p<\infty$
via the following technical fact,  \cite[Proposition 2.b.3]{LT}.

\begin{lemma}\label{l-4.1}
Let $X$ be a r.i.\ space  such that 
$0<\alpha<\underline{\alpha}_X\le \overline{\alpha}_X<\beta<1$.
Then there exist $p,q$ satisfying $1/\beta<p<q<1/\alpha$ such that
$L^q\subseteq X \subseteq L^p$ with
continuous  inclusions.
\end{lemma}

This fact yields the following equality (as linear subspaces of $L^1$):
$$
\bigcup_{1<p<\infty} L^p
=\bigcup_{0<\underline{\alpha}_X\le \overline{\alpha}_X<1} X.
$$


Recall from \eqref{3.4} that the function $w$  is
related to the description of the kernel of  $T$. 
Let $X$ be a r.i.\ space  on $(-1,1)$ with non-trivial Boyd indices.
Since $L^{2,\infty}$ is the smallest r.i.\ space  
containing $1/w$ (cf.\ Section  \ref{S2}), it follows from Lemma \ref{l-4.1}, 
that either  $T_X$ is injective or  $\mathrm{dim}(\mathrm{Ker}(T_X))=1$, depending
on whether or not $L^{2,\infty}\subseteq X$.

To indicate that $T\colon X\to X$ continuously  we simply write $T_X$, that is, $T_X\colon X\to X$.
Since $\underline{\alpha}_{X'}=1-\overline{\alpha}_X$
and $\overline{\alpha}_{X'}=1-\underline{\alpha}_X$, the condition 
$0<\underline{\alpha}_X\le \overline{\alpha}_X<1$ implies that 
$0<\underline{\alpha}_{X'} \le \overline{\alpha}_{X'}<1$;
see \cite[Proposition III.5.13]{BS}. Hence, $T_{X'}\colon X'\to X'$ is also continuous.

Recall the Parseval formula \eqref{parseval} and the Poincar\'e-Bertrand 
formula \eqref{poincare}, for a suitable pair of functions $f,g\in L^1$,
and their  importance for studying the FHT. It was noted in Section \ref{S3} that 
these formulae hold, in particular, for all pairs $f\in L^p$ and $g\in L^{p'}$ whenever
$1/p+1/p'=1$. This result was extended to any pair of functions
$f\in X$ and $g\in X'$ and all  r.i.\ spaces $X$   on $(-1,1)$
with non-trivial Boyd indices, \cite[Proposition 3.1]{COR-ampa}. 
Using the Parseval formula \eqref{parseval}  for functions $f\in X$ and $g\in X'\subseteq X^*$
it can be shown  that  the restriction
of the adjoint operator $(T_X)^*\colon X^*\to X^*$ of $T_X$ to the  associate
space $X'$ (which is a space of \textit{functions} on $(-1,1)$) is precisely $-T_{X'}\colon X'\to X'$.

The validity of both the Parseval  and the Poincar\'e-Bertrand  
formulae just mentioned have recently been  extended in the following result
to suitable pairs of functions $f\in L^1$ and $g\in\logl$,
\cite[Theorems 3.1, 3.2]{COR-asnsp}. 
For details concerning the Zygmund space $\logl$, see Section \ref{S7}.

\begin{theorem}\label{t-4.2}
Let  the functions $f\in L^1$ and $g\in \logl$ satisfy   $fT(g\chi_A)\in L^1$, 
for every  set $A\in\mathcal B$.
Then   the  Parseval formula \eqref{parseval} and the Poincar\'e-Bertrand formula 
\eqref{poincare} are valid.
\end{theorem}


The extended version of the Poincar\'e-Bertrand formula given in Theorem \ref{t-4.2}
allows the extension of  \eqref{3.4}, which  identifies $\mathrm{Ker}(T_p)\subseteq L^p$,
for $1<p<\infty$, to the operator $T$ acting in $\logl$, \cite[Theorem 3.4.]{COR-asnsp}. Namely:

\begin{theorem}\label{t-4.3}
Let $f\in \logl$. Then  $T(f)=0$ in $L^1$
if and only if $f=C/w$, for some 
constant $C\in\C$.
\end{theorem}

Regarding Theorems \ref{t-4.2} and \ref{t-4.3}, note that 
$\bigcup_{0<\underline{\alpha}_X\le \overline{\alpha}_X<1} X=
\bigcup_{1<p<\infty} L^p\subsetneqq\logl$, as shown by the function
$f(x)=(1/|x|)\log^{-\gamma}(2/|x|)$, for $x\in(-1,1)$ and any $\gamma>2$.


\section{Inversion of the finite Hilbert transform and the airfoil equation}
\label{S5}

Theorems \ref{t-3.2} and \ref{t-3.3} describe the action of the FHT on the
spaces $L^p$ for $1<p<2$ and $2<p<\infty$, respectively.
These results can be  extended to the larger classes of
r.i.\ spaces $X$ satisfying 
$1/2<\underline{\alpha}_X\le \overline{\alpha}_X<1$
and $0<\underline{\alpha}_X\le \overline{\alpha}_X<1/2$, respectively.
The main tools needed  
are various results on the continuity of the Hilbert transform in weighted $L^p$ spaces,  
\cite[Theorem 1.4.1]{GK-1}, Lemma \ref{l-4.1},   and 
Boyd's interpolation theorem, \cite[Theorem 2.b.11]{LT}. 
With these results it can be
shown that  the operator $\widehat{T}$ defined in \eqref{3.6} 
maps $X$ continuously into $X$
whenever $1/2<\underline{\alpha}_X\le \overline{\alpha}_X<1$, and that the operator
$\widecheck{T}$ defined in \eqref{3.7}  maps $X$ continuously  into $X$ 
whenever $0<\underline{\alpha}_X\le \overline{\alpha}_X<1/2$;
these operators are denoted by $\widehat{T}_X$ and $\widecheck{T}_X$,
respectively. The relevant  theorems needed for $T_X$ in this setting are the following ones, 
\cite[Theorems 3.2 and 3.3]{COR-ampa}.


\begin{theorem}\label{t-5.1}
Let $X$ be a  r.i.\ space  satisfying
$1/2<\underline{\alpha}_X\le \overline{\alpha}_X<1.$
\begin{itemize}
\item[(i)] $\mathrm{Ker}(T_X)=\mathrm{span}\{1/w\}$.
\item[(ii)] The linear operator $\widehat{T}_X$  
maps $X$ continuously into $X$ and satisfies 
 $T_X\widehat T_X=I_X$. Moreover,
$$
\int_{-1}^1\widehat{T}_X(f)(x)\,dx=0,\quad f\in X. 
$$
\item[(iii)] The operator $T_X\colon X\to X$ is surjective.
\item[(iv)] The identity $\widehat T_XT_X=I_X -P_X$ holds, 
with $P_X$ the continuous projection given by
$$
f\mapsto P_X(f):=\left(\frac1\pi\int_{-1}^1 f(t)\,dt\right)
\frac{1}{w},\quad f\in X.
$$
\item[(v)] The operator $\widehat{T}_X$ is  an isomorphism 
onto its range $\text{R}(\widehat{T}_X)$. Moreover,
$$
\text{R}(\widehat{T}_X)=\left\{f\in X: \int_{-1}^{1} f(x)\,dx=0\right\}.
$$
\item[(vi)] The following direct sum decomposition of $X$ holds: 
$$
X =\text{R}(\widehat{T}_X) \oplus \mathrm{span}\{1/w\}.
$$
\end{itemize}
\end{theorem}

Regarding the definition  of $\widecheck T$ in \eqref{3.7},   note that 
whenever  $X$ satisfies 
$0<\underline{\alpha}_X\le \overline{\alpha}_X<1/2$,
then $X'$ satisfies $1/2<\underline{\alpha}_{X'}\le \overline{\alpha}_{X'}<1$ 
and so $1/w\in X'$. Hence, for every $f\in X$, the function $f/w\in L^1$.

\begin{theorem}\label{t-5.2}
Let $X$ be a  r.i.\ space  satisfying
$0<\underline{\alpha}_X\le \overline{\alpha}_X<1/2.$
\begin{itemize}
\item[(i)] The operator $T_X\colon X\to X$ is injective.
\item[(ii)] The linear operator $\widecheck{T}_X$ 
is continuous from  $X$ into $X$ and satisfies $\widecheck T_XT_X=I_X$.
\item[(iii)] The identity $T_X\widecheck T_X=I_X -Q_X$ holds, with $Q_X$ 
the continuous projection  given by
$$
f\in X\mapsto Q_X(f):=\left(\frac1\pi\int_{-1}^1 
\frac{f(x)}{w(x)}\,dx\right) \mathbf{1}, \quad f\in X.
$$
\item[(iv)] The range of $T_X$ is the closed subspace of $X$ given by
$$
R(T_X)=\left\{f\in X: \int_{-1}^{1} \frac{f(x)}{w(x)}dx=0\right\}
=\mathrm{Ker}(Q_X).
$$
Moreover,  $\widecheck T_X$ is an isomorphism from $\text{R}(T_X)$ onto $X$.
\item[(v)] The following direct sum decomposition of $X$ holds:
$$
X= R(T_X)\oplus \mathrm{span}\{\mathbf{1}\}.
$$
\end{itemize}
\end{theorem}

Theorems \ref{t-5.1} and \ref{t-5.2} lead to the following general result on the 
inversion of the airfoil equation, \cite[Corollary 3.5]{COR-ampa}.

\begin{corollary}\label{c-5.3}
Let $X$ be a  r.i.\ space. 
\begin{itemize}
\item[(i)] Suppose that 
$1/2<\underline{\alpha}_X\le \overline{\alpha}_X<1$ and $g\in X$ is fixed. Then all 
solutions  $f\in X$ of the airfoil equation \eqref{airfoil} are given by  
$$
f=\frac{-1}{w}\; T_X (wg) + \frac{C}{w}, 
$$
with $C\in\mathbb{C}$ arbitrary. 
\item[(ii)] Suppose that  
$0<\underline{\alpha}_X\le \overline{\alpha}_X<1/2$ and $g\in X$ satisfies 
$\int_{-1}^{1} \frac{g(x)}{w(x)}dx=0.$
Then there is a unique solution $f\in X$ of the airfoil equation \eqref{airfoil}, namely  
$$
f=-w\; T_X\Big(\frac{g}{w}\Big).
$$
\end{itemize}
\end{corollary}

As a consequence of Theorems  \ref{t-5.1} and \ref{t-5.2} and Corollary \ref{c-5.3},
the operator $T_X\colon X\to X$ is well understood and there is available an inversion formula
for all r.i.\ spaces $X$ with non-trivial Boyd indices, except for  those $X$ satisfying
$\underline{\alpha}_X\le 1/2 \le \overline{\alpha}_X$. 
This class includes, for example, the Lorentz spaces $L^{2,q}$ with $1\le q\le\infty$
and, in particular, $L^2$.


\section{Extension of the finite Hilbert transform}
\label{S6}

Kolmogorov's theorem states that the FHT operator $T\colon L^1\to L^{1,\infty}$
is continuous, \cite[Theorem III.4.9(b)]{BS}. Moreover,  
$T(L^1)\not\subseteq L^1$. Hence, for any r.i.\
space $X$   necessarily  $T(L^1)\not\subseteq X$. On the other hand,
if $X$ has non-trivial Boyd indices, then
$T(X)\subseteq X$ continuously. The extension problem addresses the following question. 
Do there exist any other B.f.s.'
$Z\subseteq L^1$ such that $X\subsetneqq Z$ and $T(Z)\subseteq X$?
That is to say, given a r.i.\ space $X$ with non-trivial Boyd indices
is it possible or not to extend the finite Hilbert transform 
$T_X\colon X\to X$ to a strictly larger domain while still maintaining its values in $X$?

For $1<p<\infty$ with  $p\not=2$, Theorems \ref{t-3.2} and \ref{t-3.3} show that
$T_p\colon L^p\to L^p$ is a Fredholm operator. 
Based on this fact, it was observed in  \cite[Example 4.21]{ORS-P}, for $1<p<\infty$ 
with  $p\not=2$, that $T_p\colon L^p\to L^p$ cannot be extended to a 
strictly larger B.f.s.

As a  consequence of the inversion results Theorems \ref{t-5.1} and \ref{t-5.2}, 
the non-extendability of the FHT
was also shown to hold in those r.i.\ spaces $X$ satisfying
$1/2<\underline{\alpha}_X\le \overline{\alpha}_X<1$ and
$0<\underline{\alpha}_X\le \overline{\alpha}_X<1/2$,
\cite[Theorem 4.7]{COR-ampa}. A proof of the non-extendability of the FHT,
particular to $L^2$ (and based on its Hilbert space structure), was given in 
\cite[Theorem 5.3]{COR-ampa}.
These results left unanswered the case  when  $X$ is a r.i.\ space satisfying
$\underline{\alpha}_X\le 1/2 \le \overline{\alpha}_X$ with $X\not=L^2$.
This was settled in the following result via a unifying proof covering
all cases,  \cite[Theorem]{COR-mn}.

\begin{theorem}\label{t-6.1}
Let $X$ be any r.i.\ space on $(-1,1)$ with non-trivial Boyd indices. 
The finite Hilbert transform
$T_X\colon X\to X$ has no continuous, $X$-valued extension
to any genuinely larger B.f.s.\ containing $X$.
\end{theorem}

The proof of Theorem \ref{t-6.1} (and  of  all other results on the non-extendibility
of the FHT) relies ultimately on showing  that 
\begin{equation*}\label{5.0}
f\mapsto \sup_{|\theta|=1}\left\|T(\theta f)\right\|_X
\end{equation*}
is a norm which is equivalent to the usual norm in $X$. 
To establish this, a two-step strategy is followed. 
First, given a r.i.\ space $X$ with non-trivial Boyd indices,
a detailed study is made of the significance for a function $f\in L^1$ 
to possess the  property that  $T(f\chi_A)\in X$ for every $A\in\mathcal{B}$, 
\cite[Proposition 4.1]{COR-ampa}.

\begin{proposition}\label{p-6.2}
Let $X$ be a  r.i.\ space  with non-trivial Boyd indices. 
Given $f\in L^1$, the following conditions are equivalent.
\begin{itemize}
\item[(i)] $T(f\chi_A)\in X$ for every $A\in\mathcal{B}$.
\item[(ii)] $\displaystyle \sup_{A\in\mathcal{B}}\|T(f\chi_A)\|_X<\infty.$
\item[(iii)] $T(h)\in X$ for every $h\in L^0$ with $|h|\le |f|$ a.e.
\item[(iv)] $\displaystyle \sup_{|h|\le|f|}\|T(h)\|_X<\infty.$
\item[(v)] $T(\theta f)\in X$ for every $\theta\in L^\infty$ with $|\theta|=1$ a.e.
\item[(vi)] $\displaystyle \sup_{|\theta|=1}\|T(\theta f)\|_X<\infty.$
\end{itemize}
Moreover, if any one of $\mathrm{(i)}$-$\mathrm{(vi)}$ is satisfied, then
$$
\sup_{A\in\mathcal{B}}\big\|T(\chi_A f)\big\|_X
\le  
\sup_{|\theta|=1}\big\|T(\theta f)\big\|_X
\le 
\sup_{|h|\le|f|}\big\|T(h)\big\|_X
\le  
4 \sup_{A\in\mathcal{B}}\big\|T(\chi_A f)\big\|_X.
$$
\end{proposition}

The proof of the previous proposition given in \cite{COR-ampa}  
relies on  a deep result of Talagrand  concerning 
$L^0$-valued measures and on the Dieudonn\'e-Grothendieck Theorem for bounded 
vector measures.

The second step in the proof of Theorem \ref{t-6.1} is the 
following identification of the \textit{largest} B.f.s.\  containing $X$ 
to which $T_X\colon X\to X$ can be continuously extended, \cite[Theorem 4.6]{COR-ampa}.

\begin{proposition}\label{p-6.3}
Let $X$ be a   r.i.\ space  with non-trivial Boyd indices. 
The largest B.f.s.\ containing $X$, 
to which $T_X\colon X\to X$ can be continuously extended
while maintaining its values in $X$, is
$$
[T,X]:=\big\{f\in L^1: T(h)\in X, \;\forall |h|\le|f|\big\}
$$
equipped with the norm
\begin{equation}\label{6.1}
\|f\|_{[T,X]}:=\sup_{|h|\le|f|} \|T(h)\|_X<\infty,\quad f\in[T,X].
\end{equation}
\end{proposition}
It should be remarked that establishing the completeness and the Fatou property of 
$[T,X]$, for the norm $\|\cdot\|_{[T,X]}$, requires  some effort; see \cite[Section 4]{COR-ampa}.

Equipped with Propositions \ref{p-6.2} and \ref{p-6.3} as the main tools, the proof of the general
non-extendibility result in Theorem  \ref{t-6.1} proceeds as follows. 
For an arbitrary  $\mathcal{B}$-measurable  simple function
$$
\phi:=\sum_{n=1}^N a_n\chi_{A_n}
$$
it is clear from \eqref{6.1} that $\|\phi\|_{[T,X]}\le \|T_X\|\cdot\|\phi\|_X$. 
It is more  difficult to show that $M \|\phi\|_X \le\|\phi\|_{[T,X]}$, 
for  a constant $M>0$ depending exclusively on $X$.
To prove this, the Khintchine inequality is applied in the space
$L^1(\Lambda,\sigma)$, where 
$\Lambda:=\{1,-1\}^N$ and $\sigma$  is the product measure  of $N$ copies of the 
uniform probability measure on $\{1,-1\}$. Within this proof,
a consequence of the Stein-Weiss formula for
the distribution function of the Hilbert transform $H$ on $\R$ 
of a characteristic function (due to Laeng, 
\cite[Theorem 1.2]{La}), is crucial. Namely, for $A\subseteq\mathbb{R}$ a 
measurable set with $m(A)<\infty$,  it is the case that
$$
m\big(\{x\in A:\left| H(\chi_A)(x))\right|>\lambda\}\big) = 
\frac{2m(A)}{e^{\pi\lambda}+1},\quad \lambda>0.
$$


Combining Theorem \ref{t-6.1} and Proposition \ref{p-6.3} yields the following fact.

\begin{corollary}\label{c-6.4}
Let $X$ be a   r.i.\ space  with non-trivial Boyd indices. 
Then $X=[T,X]$ isomorphically as B.f.s.'
\end{corollary}

In the course of the  above investigations the following (rather unexpected) 
characterization of when a function $f\in L^1$  belongs to
$X$, in terms of the set of its $T$-transforms $\{T(f\chi_A): A\in\mathcal{B}\}$, 
was established, \cite[Corollary]{COR-mn}.

\begin{proposition}\label{p-6.5}
Let $X$ be a r.i.\ space on $(-1,1)$ with non-trivial Boyd indices. 
\begin{itemize}
\item[(i)] 
For a function $f\in L^1$ the following conditions are equivalent.
\begin{itemize}
\item[(a)] $f\in X$.
\item[(b)] $T(f\chi_A)\in X$ for every $A\in\mathcal{B}$.
\item[(c)] $T(f\theta)\in X$ for every  $\theta\in L^\infty$ with $|\theta|=1$ a.e.
\item[(d)] $T(h)\in X$ for every $h\in L^0$ with $|h|\le |f|$ a.e.
\end{itemize}
\item[(ii)]
There exists a constant $\beta>0$ such that, for every $f\in X$, we have 
$$
\frac{\beta}{4}\|f\|_X
\le  \sup_{A\in\mathcal{B}}\big\|T(\chi_A f)\big\|_X
\le \sup_{|\theta|=1}\big\|T(\theta f)\big\|_X
\le \sup_{|h|\le|f|}\big\|T(h)\big\|_X
\le \|T_X\|\cdot\|f\|_X .
$$
\end{itemize}
\end{proposition}


\section{The finite Hilbert transform acting on $\logl$}
\label{S7}

In all investigations so far $T$ was always considered as a linear operator acting
from a r.i.\  space  into \textit{itself}. We now consider $T$ when it is acting in the classical
Zygmund space  $\logl:=\logl(-1,1)$. As will become clear, $T(\logl)\not\subseteq \logl$.

The Zygmund space $\logl$  consists of all measurable functions $f$
on $(-1,1)$ for which either one of the following two equivalent conditions hold:
$$
\int_{-1}^1|f(x)|\log^+|f(x)|\,dx<\infty,\quad \int_0^2 f^*(t)\log\Big(\frac{2e}{t}\Big)\,dt<\infty;
$$
see \cite[Definition IV.6.1 and Lemma IV.6.2]{BS}. 
The space $\logl$ is  r.i.\  with a.c.-norm (cf. \cite[p.247-248]{BS})  given by
$$
\|f\|_{\logl}:=\int_0^2 f^*(t)\log\Big(\frac{2e}{t}\Big)\,dt,\quad f\in\logl.
$$
Then $\logl$ is a r.i.\  space on $(-1,1)$ close to $L^1$ in the sense that
$L^p\subseteq\logl$ for all $1<p<\infty$, \cite[Theorem IV.6.5]{BS}, 
which implies (in view of Lemma \ref{l-4.1}) that 
$X\subseteq \logl$ for all r.i.\   spaces  $X$  with non-trivial Boyd indices. 
The associate space of $\logl$ is the space $\lexp$ consisting 
of all measurable functions $f$ on $(-1,1)$ having exponential integrability;
see \cite[Definition IV.6.1 \& Theorem IV.6.5]{BS}. 
The separability of  $\logl$ implies that $(\logl)^*=(\logl)'=\lexp$.

The Boyd indices of $\logl$ are \textit{trivial}, namely,
$\underline{\alpha}_{\logl}= \overline{\alpha}_{\logl}=1$,
\cite[Theorem IV.6.5]{BS}, and so $T$ \textit{cannot} map $\logl$ into itself.
However, it turns out that
$T\colon\logl\to L^1$ is a continuous operator, \cite[Theorem 2.1]{COR-asnsp}.
Whenever convenient, the operator $T\colon\logl\to L^1$ will
also be denoted by $\tlog$.

Note that Theorems \ref{t-4.2} and \ref{t-4.3} imply the next result.

\begin{corollary}\label{c-7.1}
The  Parseval  formula 
$$
\int_{-1}^1fT(g)=-\int_{-1}^1gT(f),\quad f\in L^\infty,\; g\in\logl, 
$$
holds, as does the Poincar\'e-Bertrand formula (pointwise in $L^0$) 
$$
T(gT(f)+fT(g))=(T(f))(T(g))-fg,\quad f\in L^\infty,\; g\in\logl.
$$
Moreover,  $\mathrm{Ker}(\tlog)=\mathrm{span}\{1/w\}$.
\end{corollary}


In order to study the inversion of the FHT
on $\logl$,   the operator $\widehat T$ defined in \eqref{3.6} 
plays a central role (as it did in Theorems \ref{t-3.2} and \ref{t-5.1} for the 
inversion of $T$ when acting on the $L^p$-spaces, for $1<p<2$, and on r.i.\ spaces 
$X$ satisfying $1/2<\underline{\alpha}_X\le \overline{\alpha}_X<1$, respectively). 
This is also the case when $T$ acts in the space $\logla$, for each $\alpha>1$,
consisting of all measurable functions $f$
on $(-1,1)$ for which either one of the following two equivalent conditions hold:
$$
\int_{-1}^1|f(x)|(\log(2+|f(x)|)^\alpha\,dx<\infty,\quad \int_0^2 f^*(t)\log^\alpha
\Big(\frac{2e}{t}\Big)\,dt<\infty.
$$
The space $\logla$ is  r.i.\  for the a.c.-norm given by
$$
\|f\|_{\logla}:=\int_0^2 f^*(t)\log^\alpha\Big(\frac{2e}{t}\Big)\,dt,\quad f\in\logla;
$$
see \cite[Definition IV.6.11 and Lemma IV.6.12]{BS}.  
The following inclusions hold: 
$$
L^p\subseteq \loglb\subseteq \logla\subseteq \logl,
\quad 1<p,\,\, 1<\alpha<\beta.
$$

A new feature is that extrapolation results   enter in the study of $\tlog$.
In particular, the following theorem of Yano
is important, \cite[Theorem XII.(4.41)]{Z}.

\begin{theorem}\label{t-7.2}
Let $1<p_0<\infty$ and  $S$ be a linear  operator that 
maps $L^p$ continuously into $L^p$
for all $1<p<p_0$ and  such that there exist  constants $C>0$ and 
$1<\alpha\le p_0$ satisfying 
$$
\|S\|_{L^p\to L^p}\le \frac{C}{p-1},\quad p\in(1,\alpha).
$$
Then $S$ can be extended  to $\logl$ with $S\colon \logl\to L^1$
a continuous operator.
\end{theorem}

The following related result,  \cite[Theorem 5.1]{EK}, is also required.

\begin{theorem}\label{t-7.3}
Let $p_0, S, C$ and $\alpha$ satisfy the conditions in Yano's theorem.
Suppose, for some $\gamma\ge0$, that 
$S\colon L(\textnormal{log } L)^{\gamma}\to L^1$ 
is continuous. Then,  for every $\beta>0$, also 
$$
S\colon L(\textnormal{log } L)^{\gamma+\beta}\to L(\textnormal{log } L)^{\beta}   
$$
continuously.
\end{theorem}

Pichorides'  calculation
of the norm of the Hilbert transform $H\colon L^p(\R)\to L^p(\R)$,
\cite{P}, together with a result of McLean and Elliott showing that
$\|T\|_{L^p\to L^p}=\|H\|_{L^p(\R)\to L^p(\R)}$, \cite[Theorem 3.4]{ME}, imply  that 
$$
\|T\|_{L^p\to L^p}  = \tan(\pi/(2p)) \le\frac{3}{p-1},\quad 1<p<2;
$$
see \cite[Lemma 4.3]{COR-asnsp}.
Combined with Theorem \ref{t-7.3}, this yields the following result,
\cite[Proposition 4.4]{COR-asnsp}.

\begin{proposition}\label{p-7.4}
The finite Hilbert transform $T$ satisfies
\begin{equation*}
T\colon L(\textnormal{log } L)^{1+\beta}\to L(\textnormal{log } L)^{\beta}
\end{equation*}
continuously, for every  $\beta\ge0$.
\end{proposition}

The proof of the corresponding result for the operator
$\widehat T$ defined in \eqref{3.6}   follows a similar strategy. However,
it requires not only  that  $\widehat{T}\colon L^p\to L^p$ continuously, for $1<p<2$, 
(a special case of Khvedelidze's Theorem  \ref{t-3.1}), 
 but also an \textit{explicit} upper bound on the operator norms
$\|\widehat{T}\|_{L^p\to L^p}$ for $p$ near 1, which is not given in \cite{GK-1}, \cite{MP}.
This upper bound is established 
via some technical auxiliary facts, \cite[Lemmas 4.6, 4.7, 4.8 and 6.1]{COR-asnsp},
which then lead to a proof of the following result, \cite[Proposition 4.5]{COR-asnsp}.


\begin{proposition}\label{p-7.5}
For each $\beta\ge0$ the operator
\begin{equation*}
\widehat{T}\colon L(\textnormal{log } L)^{1+\beta}\to L(\textnormal{log } L)^{\beta}
\end{equation*}
continuously. In particular, $\widehat{T}\colon \logl\to L^1$ is continuous.
\end{proposition}


For $\tlog\colon\logl\to L^1$, the result corresponding to 
Theorems \ref{t-3.2} and \ref{t-3.3} for the $L^p$-spaces and to 
Theorems \ref{t-5.1} and \ref{t-5.2} for  r.i.\ spaces 
is the following one, \cite[Theorem 4.10]{COR-asnsp}.

\begin{theorem}\label{t-7.6}  
The following assertions are valid.
\begin{itemize}
\item[(i)]   The operator $\tlog\colon\logl\to L^1$ is not injective.
\item[(ii)] Let  $g\in\logl$.   Then $\widehat T(g)\in L^1$ and $T(\widehat T(g))=g$. Moreover,
\begin{equation*}
\int_{-1}^1\widehat{T}(g)(x)\,dx=0.
\end{equation*}
\item[(iii)] The operator $P\colon\logl\to\logl $  given by
\begin{equation*}
P (f)(x):=\left(\frac1\pi\int_{-1}^1 f(t)\,dt\right)
\frac{1}{w},   \quad f\in {\logl} ,
\end{equation*}
is a continuous projection satisfying the inequality
$$
\|P\|_{\logl\to\logl}\le \frac1\pi\Big\|\frac{1}{w}\Big\|_{\logl}.
$$
Furthermore, $TP\colon\logl\to L^1$ is the zero operator.
\item[(iv)] For each $f\in\logl$ it is the case that 
$$
f-P(f)=\widehat{T}(T(f)).
$$
Moreover, $\widehat{T}T\colon\logl\to\logl$ satisfies the inequality
$$
\|\widehat{T}T\|_{\logl\to\logl}
\le 1+\frac1\pi\Big\|\frac{1}{w}\Big\|_{\logl}.
$$
\end{itemize}
\end{theorem}


Theorem \ref{t-7.6} yields a description of the range space of $\tlog$,
\cite[Corollary 4.11]{COR-asnsp}.

\begin{corollary}\label{c-7.7}  
A function $g\in L^1$ belongs to the range space
$T(\logl)$ of $\tlog$  if and only if it satisfies  both $\widehat T(g)\in \logl$ and
$T(\widehat{T}(g))=g$. That is, 
$$
T(\logl)=\Big\{g\in L^1:  \widehat T(g)\in \logl,\; T(\widehat{T}(g))=g\Big\}.
$$
\end{corollary}

This description of the range space $T(\logl)$ is precise and rather useful, despite
not being  fully explicit. Additional properties of  $T(\logl)$ are presented in the next result,
\cite[Proposition 4.12]{COR-asnsp}.

\begin{proposition}\label{p-7.8} 
The following assertions hold for the continuous linear operators  $\tlog\colon \logl \to L^1$
and $\widehat T\colon \logl \to L^1$.
\begin{itemize}
\item[(i)] The range  $T(\logl)$ of $\tlog$ is a proper, dense, linear subspace of $L^1$.
\item[(ii)]  $\logld$ is included in the range $T(\logl)$.
\item[(iii)] Neither $T(\logl)$ nor $\widehat T(\logl)$ is contained in $\logl$.
\item[(iv)]   $\logl$ is not included in the range $T(\logl)$  
  nor in the range  $\widehat T(\logl)$. 
\end{itemize}
\end{proposition}

The information available in Theorem \ref{t-7.6} and Corollary \ref{c-7.7}
leads to an inversion formula for the operator $\tlog\colon\logl\to L^1$, 
\cite[Theorem 4.14]{COR-asnsp}.

\begin{theorem}\label{t-7.9}
Let $g$ belong to the range of  $\tlog\colon\logl\to L^1$.  
All solutions  $f\in \logl$ of the airfoil equation \eqref{airfoil} are of the form   
$$
f=\frac{-1}{w}\; T(w g) + \frac{C}{w}, 
$$
with  $C\in\mathbb{C}$ arbitrary.  In particular,
this applies to   every $g\in \logld\subseteq T(\logl)$.
\end{theorem}


As was done for $T_X\colon X\to X$, we now discuss the extension problem for the
operator $\tlog\colon \logl \to L^1$. A similar approach as in Section \ref{S6} is applicable.
That is,  it is again possible to identify  the \textit{largest} B.f.s.\  containing $\logl$ 
to which $\tlog\colon \logl \to L^1$ can be continuously extended. Namely, 
$$
[T, L^1] : = \Big\{f \in L^1: T(h) \in L^1 \text{ for all } |h|  \le |f|\Big\},
$$
together with its associated functional
$$
\|f\|_{[T, L^1] }:= \sup\Big\{\|T(h)\|_{L^1}: |h|\le|f| \Big\},\quad f \in [T, L^1] ,
$$
\cite[Lemmas 5.4 and 5.5]{COR-asnsp}.
To verify the inclusion $\logl \subseteq [T, L^1] $ is direct and so
it remains to establish the opposite containment  $[T, L^1]\subseteq  \logl$.
This is done in \cite[Theorem 5.3]{COR-asnsp}, and is based on  
a result of Stein concerning the space $\logl$, \cite[Theorem 3(b)]{S}.
Combining these facts yields the following theorem, \cite[Theorem 5.6]{COR-asnsp}.

\begin{theorem}\label{t-7.10}
The identity $[T, L^1] = \logl$ holds as an order 
and bicontinuous isomorphism between B.f.s.' Consequently, 
$\tlog\colon\logl \to L^1$ does not admit a continuous 
$L^1$-valued linear  extension to 
any strictly  larger B.f.s.\  within $ L^1$ and containing $\logl$.
\end{theorem}

It should be noted, however, that there do exist functions in
$L^1\setminus \logl$  which $T$ maps  into $\logl\setminus T(\logl)$.
That is, the   linear subspace 
$T^{-1}(L^1)=\{f\in L^1:T(f)\in L^1\}$ of $L^1$ is strictly larger than the optimal 
domain $[T, L^1]=\logl$.
The  difference between those two  spaces 
is that  the optimal domain $[T, L^1]$  is a function lattice
(i.e., it satisfies the \textit{ideal property} 
namely, $f\in [T, L^1]$ and $|g|\le|f|$ a.e.\ imply  
that $g\in [T, L^1]$), whereas  $T^{-1}(L^1)$ is not;
see  \cite[Proposition 5.8 \& Remark 5.9]{COR-asnsp}.


\section{The fine spectrum of the finite Hilbert transform}
\label{S8}

The aim of this section is to give a detailed expose 
of the fine spectra of the FHT 
acting in r.i.\ spaces over $(-1,1)$ with non-trivial Boyd indices 
which goes  beyond that known to date for the classical $L^p$-theory.

The  spectrum of $T_p\colon L^p\to L^p$, for $1<p<\infty$, was completely identified
by Widom in 1960, who also described  its fine spectra, that is, the point spectrum, 
continuous spectrum and residual spectrum, \cite{W}; see
also \cite[\S13.6]{J}. It is worthwhile to describe Widom's results, for which   
we will use the following minor modification of the FHT (\textit{only} in this section), namely,
\begin{equation*}\label{T2}
T(f)(t)=\lim_{\varepsilon\to0^+} \frac{1}{\pi i}
\left(\int_{-1}^{t-\varepsilon}+\int_{t+\varepsilon}^1\right) \frac{f(x)}{x-t}\,dx, 
\end{equation*}
which differs from \eqref{T} only by a factor of $1/i$.

For $1<p<\infty$, consider the subset of $\C$ given by
$$
\mathcal{R}_p:=\big\{\pm1\big\}\cup
\left\{\lambda\in\C: \frac{1}{2\pi} 
\left|\arg\bigg(\frac{1+\lambda}{1-\lambda}\bigg)\right|\le \left|\frac12-\frac1p\right|\right\},
$$
which is the region bounded by both, the circular arc with end-points
$\pm1$ which passes through $i\cot(\pi/p)$, together with the 
circular arc having end points $\pm1$ which passes through $i\cot(\pi/p')$,
where $1/p+1/p'=1$; see the following diagram. 

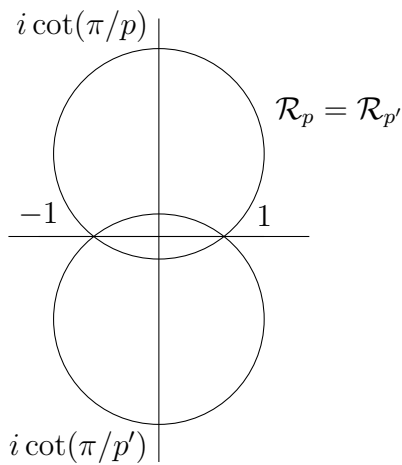
\begin{center}
\begin{tikzpicture}
\draw[-,xshift=-0cm] (-2,0) -- coordinate (x axis mid) (2,0);
\draw[-,xshift=-0cm] (0,-3) -- coordinate (y axis mid)(0,2.9);
\draw (0,-1.1) circle (1.4);
\draw (0,1.1) circle (1.4);
\draw (-1.6,0) node[above]{$-1$} -- (1.4,0)  node[above]{$1$};
\draw (0,2.8) node[left]{$i\cot(\pi/p)$};     
\draw (0,-2.8) node[left]{$i\cot(\pi/p')$};  
\draw (2.4,2) node[below]{$\mathcal{R}_p=\mathcal{R}_{p'}$} ;
\end{tikzpicture}
\end{center}

It is an important feature that $\mathcal{R}_p=\mathcal{R}_{p'}$.
Note that $\mathcal{R}_2=[-1,1]$ and, for $1<p<\infty$, that the set 
$\mathcal{R}_p$ increases as $|p-2|$ increases. 
The geometric symmetries of $\sigma(T_p)$, which can be gleaned from the
diagram above and are formulated in Proposition \ref{p-8.3} below,
will play an important role.

The following result is due to Widom; see Remark 1(2) and 
Remark 2 (pp. 156-157) in \cite{W}. The interior 
of a set $B\subseteq \C$ is denoted by $\mathrm{int}(B)$
and its boundary by $\partial B$.

\begin{theorem}\label{t-8.1}
Let $1<p<\infty$. The operator $T_p\colon L^p\to L^p$ has spectrum
$$
\sigma(T_p)=\mathcal{R}_p.
$$ 
Regarding the fine spectra of $T_{p}$ the following identifications hold.
\begin{itemize}
\item[(i)] Let $1<p<2$. Then $\sigma_{\mathrm{pt}}(T_p) =\mathrm{int}(\mathcal{R}_p)$,
$\sigma_{\mathrm{r}}(T_p)=\emptyset$ and $\sigma_{\mathrm{c}}(T_p)=\partial\mathcal{R}_p$.
\item[(ii)] Let $p=2$. Then $\sigma_{\mathrm{pt}}(T_2) =\emptyset$,
$\sigma_{\mathrm{r}}(T_2)=\emptyset$ and 
$\sigma_{\mathrm{c}}(T_2)=\partial\mathcal{R}_2=\mathcal{R}_2$.
\item[(iii)] Let $2<p<\infty$. Then $\sigma_{\mathrm{pt}}(T_p) =\emptyset$,
$\sigma_{\mathrm{r}}(T_p)=\mathrm{int}(\mathcal{R}_p)$ 
and $\sigma_{\mathrm{c}}(T_p)=\partial\mathcal{R}_p$.
\end{itemize}
\end{theorem}

As a consequence of Theorem \ref{t-8.1}, the set $\mathcal{A}$ of
all complex numbers which occur as an eigenvalue of $T$, when $T$ 
acts in $L^p$ for some $1<p<\infty$, is given by 
\begin{equation}\label{8.1}
\mathcal{A}=\bigcup_{1<p<\infty} \mathrm{int}(\mathcal{R}_p)
=\C\setminus\big\{(-\infty,-1]\cup[1,\infty)\big\}.
\end{equation}
The next result, due to J\"orgens in 1970, \cite[Theorem 13.9]{J},
identifies the set of all possible eigenfunctions for $T$ when $T$ 
acts over  $\bigcup_{1<p<\infty}L^p$.

\begin{theorem}\label{t-8.2}
For each  $\lambda\in\mathcal{A}$, the corresponding eigenspace of 
$T$ is the one-dimensional subspace $\mathrm{span}\{\xi_\lambda\}\subseteq L^p$
spanned by the eigenfunction
\begin{equation}\label{8.2}
\xi_\lambda(x):=\frac{1}{w(x)}\left(\frac{1-x}{1+x}\right)^{z(\lambda)},\quad  x\in(-1,1),
\end{equation}
for all $1<p<\infty$ such that $\xi_\lambda\in L^p$,  where the function $z(\cdot)$ is given by
\begin{equation}\label{8.3}
z(\lambda):=\frac{1}{2\pi i}\log\left(\frac{1+\lambda}{1-\lambda}\right),\qquad z(0)=0.
\end{equation}
\end{theorem}

So, the  $L^p$-theory concerning the spectrum of $T$ was completely 
determined by the 1970s. According to   Lemma \ref{l-4.1},
the previous result implies that
\begin{equation*}
\mathcal{A}=\Big\{\lambda\in\C: T(f)=\lambda f \text{ for some }   f\in X\setminus\{0\} 
\text{ and an $X$ with } 
0<\underline{\alpha}_X\le \overline{\alpha}_X<1 
\Big\},
\end{equation*}
and, with $\mathcal{E}=\{\xi_\lambda:\lambda\in\mathcal{A}\}$,  that
\begin{equation*}
\mathcal{E}=\Big\{\xi\in X: 0<\underline{\alpha}_X\le \overline{\alpha}_X<1,  
\;\exists \lambda\in\mathcal{A} \text{ such that $\xi=\xi_\lambda$}
\,(\textrm{cf.}  \eqref{8.2})\Big\}.
\end{equation*}

A subset $D\subseteq\C$ is called \textit{$\R$-balanced}
if $\alpha\lambda\in D$ whenever $\lambda\in D$ and
$\alpha\in\R$ satisfies $|\alpha|\le1$. The following result
(a combination of Proposition 3.1 and Corollaries 5.3, 5.4 and 5.6 in
\cite{COR-am}) indicates the strong symmetry properties of the fine spectra of $T$.

\begin{proposition}\label{p-8.3}
Let $X$ be any   r.i.\ space  on $(-1,1)$  
with non-trivial Boyd indices. 
\begin{itemize}
\item [(i)] Each of the spectra $\sigma(T_X)$, 
$\sigma_{\mathrm{pt}}(T_X)$, $\sigma_{\mathrm{c}}(T_X)$ and 
$\sigma_{\mathrm{r}}(T_X)$  is symmetric with respect to both 
reflection in the real axis and in the imaginary axis in $\C$.
In particular, these spectra are also symmetric with respect to reflection through 0.
\item [(ii)] The set $\sigma_{\mathrm{pt}}(T_{X})$ is $\R$-balanced.
In addition, if $X$ is separable, then also 
$\sigma_{\mathrm{r}}(T_{X})$ is $\R$-balanced.
\item [(iii)] It is always the case that 
$[-1,1]\subseteq \sigma(T_{X})$. If, in addition, $X$ is separable, then
$\pm1\in\sigma_{\mathrm{c}}(T_{X})$.
\end{itemize}
\end{proposition}

The proof of parts (i) and (iii) in Proposition \ref{p-8.3} is essentially via
manipulations of the definitions involved, whereas part (ii) makes
explicit use of properties of the functions $\xi_\lambda$ that occur in
\eqref{8.2}. The proof that $\sigma_{\mathrm{pt}}(T_{X})$
is $\R$-balanced  requires an analysis of the decreasing 
rearrangement $\xi_\lambda^*$ of the
eigenfunctions  $\xi_\lambda$ in \eqref{8.2}.

For the spectrum of a (continuous)  Banach space operator $A$ it is known that 
$$
\sigma_{\mathrm{r}}(A)\subseteq \sigma_{\mathrm{pt}}(A^*),\quad 
\sigma_{\mathrm{pt}}(A^*) \subseteq\sigma_{\mathrm{pt}}(A)\cup
\sigma_{\mathrm{r}}(A),\quad
\sigma_{\mathrm{pt}}(A)\subseteq \sigma_{\mathrm{pt}}(A^*)
\cup\sigma_{\mathrm{r}}(A^*),
$$
\cite[Theorem 5.13]{J}. Recall if $X$ is a separable r.i.\ space on $(-1,1)$, then
its associate space $X'$ equals $X^*$ with 
$\underline{\alpha}_{X'}=1-\overline{\alpha}_{X}$
and $\overline{\alpha}_{X'}=1-\underline{\alpha}_{X}$; see
Section \ref{S4}. Moreover, whenever 
$0<\underline{\alpha}_{X}\le\overline{\alpha}_{X}<1$ 
we recall that the restriction of $(T_X)^*$ to $X'$ is precisely $-T_{X'}$. 
These observations indicate the usefulness
of duality arguments, which are often used in \cite{COR-am},
when they are combined with the symmetry properties in Proposition \ref{p-8.3}.

Proposition \ref{p-8.3}(iii) shows that \textit{always} 
$0\in\sigma(T_X)$. The following result,
\cite[Proposition 5.1]{COR-am}, indicates the distinguished role played by the location 
of the point $0$ within $\sigma(T_X)$. Recall 
the special relevance of the spaces $L^{2,\infty}$ and $L^{2,1}$
and the fact that $|x|^{-1/2}\in X$ if and only if $L^{2,\infty}\subseteq X$;
see Section \ref{S4}. Again the functions $\xi_\lambda$ play an
important role in the proof.

\begin{proposition}\label{p-8.4}
Let $X$ be a  separable r.i.\ space  on $(-1,1)$ with non-trivial Boyd indices. 
Precisely one of the following three mutually exclusive alternatives holds.
\begin{itemize}
\item[(i)]  The following conditions are equivalent. 
\begin{itemize}
\item[(a)] $0\in\sigma_{\mathrm{pt}}(T_X)$.
\item[(b)] $\sigma_{\mathrm{pt}}(T_X)\not=\emptyset$.
\item[(c)] $|x|^{-1/2}\in X$.
\item[(d)] $L^{2,\infty}\subseteq X$.
\item[(e)] $(-1,1)\subseteq\sigma_{\mathrm{pt}}(T_X)$.
\end{itemize}

\item[(ii)] The following conditions are equivalent. 
\begin{itemize}
\item[(a)] $0\in\sigma_{\mathrm{r}}(T_X)$.
\item[(b)] $\sigma_{\mathrm{r}}(T_X)\not=\emptyset$.
\item[(c)] $|x|^{-1/2}\in X'$.
\item[(d)] $X\subseteq L^{2,1}$.
\item[(e)] $(-1,1)\subseteq\sigma_{\mathrm{r}}(T_X)$.
\end{itemize}

\item[(iii)]  The following conditions are equivalent. 
\begin{itemize}
\item[(a)] $0\in\sigma_{\mathrm{c}}(T_X)$.
\item[(b)] $|x|^{-1/2}$ belongs to neither $X$ nor to $X'$.
\item[(c)] $\sigma(T_X)=\sigma_{\mathrm{c}}(T_X)\supseteq[-1,1]$.
\end{itemize}
\end{itemize}
\end{proposition}

The previous result, together with Remark 5.2  in \cite{COR-am},
imply the following somewhat unexpected trichotomy
(for \textit{all}  separable r.i.\ spaces $X$ with non-trivial Boyd indices):
$$
\sigma(T_X)=\sigma_{\mathrm{pt}}(T_X)\cup\sigma_{\mathrm{c}}(T_X);
\quad
\sigma(T_X)=\sigma_{\mathrm{r}}(T_X)\cup\sigma_{\mathrm{c}}(T_X);
\quad
\sigma(T_X)=\sigma_{\mathrm{c}}(T_X).
$$

The M\"obius transformation $u(\lambda):=(1+\lambda)/(1-\lambda)$
maps the set $\mathcal{A}$ of all eigenvalues given in \eqref{8.1} 
onto the set $\Omega:=\C\setminus(-\infty,0]$. In $\Omega$ the branch of the 
argument used for complex numbers is fixed to lie in $(-\pi,\pi)$.
Then, for the function $z(\cdot)$ in \eqref{8.3}, its real part is given by 
\begin{equation}\label{8.4}
\Re(z(\lambda))=\frac{1}{2\pi}\arg\Big(\frac{1+\lambda}{1-\lambda}\Big)=
\frac{1}{2\pi}\arg(u(\lambda))\in\Big(-\frac12,\frac12\Big).
\end{equation}
Given any r.i.\ space $X$ on $(-1,1)$ with non-trivial 
Boyd indices, define $p_X\in (1,\infty)$ by
\begin{equation}\label{8.5}
p_X:=\inf\Big\{p\in(1,\infty): |x|^{-1/p}\in X\Big\}=\inf\Big\{p\in(1,\infty): 
L^{p,\infty}\subseteq X\Big\},
\end{equation}
where the fact is used that $|x|^{-1/p}\in X$ if and only if
$L^{p,\infty}\subseteq X$. The index $p_X$ can be attained or not,
depending on $X$, which will be important for certain properties
of $\sigma(T_X)$. For instance, if $X=L^{p,r}$ with $1<p<\infty$
and $1\le r<\infty$, then $p_X=p$ and the infimum in \eqref{8.5} 
is not attained, whereas for $X=L^{p,\infty}$ we have that  
$p_X=p$ and the infimum in \eqref{8.5} is  attained.  
For each $\lambda\in\mathcal{A}$, define $\gamma_\lambda>0$ by
$$
\frac{1}{\gamma_\lambda} := \frac12+\frac{1}{2\pi}
\left|\arg\bigg(\frac{\lambda+1}{\lambda-1}\bigg)\right|.
$$
From the complex argument specified in \eqref{8.4} it follows that
$1<\gamma_\lambda\le2$. The following technical result, \cite[Lemma 3.2]{COR-am},
determines precisely when an eigenvector $\xi_\lambda$ belongs to $X$ in 
terms of an inequality between $\gamma_\lambda$ and the index $p_X$.

\begin{lemma}\label{l-8.5}
Let $X$ be any   r.i.\ space  on $(-1,1)$  with non-trivial Boyd indices.
Let $\lambda\in\mathcal{A}$ and $\xi_\lambda\in\mathcal{E}$ be the 
corresponding  eigenfunction. 
\begin{itemize}
\item[(i)] If $p_X$ is  attained, then 
$\xi_\lambda \in X$ if and only if $p_X\le \gamma_\lambda.$
\item[(ii)] If $p_X$ is not attained, then
$\xi_\lambda \in X$ if and only if $p_X<\gamma_\lambda.$
\end{itemize}
\end{lemma}

From Theorems \ref{t-8.1} and \ref{t-8.2} and Lemma \ref{l-8.5}
we note that $T(f)=\lambda f$, for $f\in X\setminus\{0\}$ and $\lambda\in\C$, precisely when 
$\lambda\in\mathcal{A}$ and
$f$ belongs to the one-dimensional eigenspace $\mathrm{span}\{\xi_\lambda\}$ 
spanned by $\xi_\lambda$. So, $\sigma_{\mathrm{pt}}(T_X)$
consists of those $\lambda\in\C$ for which $\xi_\lambda\in X$, that is,
$$
\sigma_{\mathrm{pt}}(T_X)=\{\lambda\in\C:\xi_\lambda\in
\mathcal{E}\cap \ X\}.
$$
This observation, together with Lemma \ref{l-8.5}, are the 
essential ingredients in the proof of the following 
characterization of the point spectrum of $T_X$, 
\cite[Proposition 3.3]{COR-am}.

\begin{proposition}\label{p-8.6}
Let $X$ be any   r.i.\ space  on $(-1,1)$  with non-trivial Boyd indices.
\begin{itemize}
\item [(i)] Let $p_X>2$ (attained or not) or,  $p_X=2$ with $p_X$  not attained.
Then $\sigma_{\mathrm{pt}}(T_X)=\emptyset$.
\item [(ii)]   Let $p_X\le2$ with $p_X$  attained. Then 
$\sigma_{\mathrm{pt}}(T_X)=\mathcal{R}_{p_X}\setminus\{\pm1\}$.
\item [(iii)] Let $p_X<2$ with $p_X$  not attained. Then
$\sigma_{\mathrm{pt}}(T_X)=\mathrm{int}(\mathcal{R}_{p_X})$.
\end{itemize}
\end{proposition}


Using the results recorded so far, together with Boyd's interpolation 
theorem, \cite[Theorem 2.b.11]{LT}, in a series of five
propositions/theorems in Section 4 of \cite{COR-am}, a full description
of the fine spectra of $T$ acting in the Lorentz spaces $L^{p,r}$,
for $1<p<\infty$ and $1\le r<\infty$, is presented. 
It should be noted  that the symmetry properties in 
Proposition \ref{p-8.3} play an important role in this description. The following table, 
where $T_{p,r}$ denotes $T\colon L^{p,r}\to L^{p,r}$, gives a
complete summary of these results; see p.16 in \cite{COR-am}.
Observe, for $p=r$, that we recover Theorem \ref{t-8.1} of Widom.

\bigskip
\begin{center}
\begin{tabular}{ |c | c| c | c |  c| }
\hline
   	$L^{p,r}$ &$\sigma(T_{p,r})=\mathcal{R}_p$
	&$\sigma_{\mathrm{pt}}(T_{p,r})$& $\sigma_{\mathrm{r}}		
	(T_{p,r})$&$\sigma_{\mathrm{c}} (T_{p,r})$
\\ \hline 
 	$1<p<2$&$1\le r<\infty$  &$\mathrm{int}(\mathcal{R}_p)$ 
	& $\emptyset$&$\partial \mathcal{R}_p$
\\ \hline
  	$2<p<\infty$& $r=1$ &$\emptyset$&$\mathcal{R}_p\setminus\{\pm1\}$&$\{\pm1\}$
\\ \hline
 	&$1< r<\infty$  & $\emptyset$&$\mathrm{int}(\mathcal{R}_p)$  
	&$\partial \mathcal{R}_p$
\\ \hline
   	$p=2$&$r=1$ &$\emptyset$&$(-1,1)$&$\{\pm1\}$
\\ \hline
 	&$1< r<\infty$  & $\emptyset$&$\emptyset$  &$[-1,1]$
\\ \hline
\end{tabular}
\end{center}
\bigskip

Proposition \ref{p-8.4} suggests a strategy of how  to investigate $\sigma(T_X)$ 
further. The key point is to decide in which part of the spectrum of $T_X$
the point 0 lies. This requires introducing an index additional to $p_X$ (cf.\ \eqref{8.5}).

Let $X$ be a   r.i.\ space  on $(-1,1)$  with non-trivial Boyd indices. 
Define  $q_X\in(1,\infty)$ by
$$
q_X:=\sup\Big\{q\in(1,\infty): X\subseteq L^{q,1}\Big\}.
$$
The index  $q_X$ can be  attained or not, depending on the space $X$. 
It is clear that $q_X\le p_X$. Note that  
there is no r.i.\ space $X$ on $(-1,1)$ with non-trivial Boyd indices
for which $p_X=q_X$ with both $p_X$ and $q_X$ being attained.
The following result (cf.\ \cite[Lemma 6.1]{COR-am}) 
presents some of the connections between the two indices $p_X$ and $q_X$. Recall from
Section \ref{S4} that $\underline{\beta}_X$ and $\overline{\beta}_X$
are the lower and   upper fundamental indices of $X$, respectively.

\begin{lemma}\label{l-8.8} 
Let $X$ be any r.i.\ on $(-1,1)$ with non-trivial Boyd indices.
The following inequalities hold:
$$
0<\underline{\alpha}_X\le \underline{\beta}_X\le
1/p_X\le 1/q_X\le 
\overline{\beta}_X\le \overline{\alpha}_X<1.
$$
For the associate space $X'$  of $X$, it is the case  that $p_{X'}=(q_X)'$. 
Moreover, $p_{X'}$ is attained if and only if $q_X$ is attained.
\end{lemma}

The following result, whose proof relies on Lemma \ref{l-8.8}, treats the cases
when $0\in\sigma_{\mathrm{pt}}(T_X)$, resp.\ $0\in\sigma_{\mathrm{r}}(T_X)$,
resp.\ $0\in\sigma_{\mathrm{c}}(T_X)$ and consists successively of Proposition
6.2, 6.3 and 6.4 in \cite{COR-am}. It provides a 
comprehensive description of the fine spectra for a large class of r.i.\ spaces.

\begin{proposition}\label{p-8.9}
Let $X$ be a separable r.i.\ space on $(-1,1)$ with non-trivial Boyd indices. 
\begin{itemize}
\item [(i)] If both $p_X<2$ and $p_X$ is not attained, then 
$$
\sigma_{\mathrm{pt}}(T_X) =\mathrm{int}(\mathcal{R}_{p_X});\;
\sigma_{\mathrm{r}}(T_X)=\emptyset,
$$
whereas if both $p_X\le2$ and $p_X$ is attained, then 
$$
\sigma_{\mathrm{pt}}(T_X)=\mathcal{R}_{p_X}\setminus\{\pm1\};\;
\sigma_{\mathrm{r}}(T_X)=\emptyset.
$$
\item [(ii)] If both  $q_X>2$ and $q_X$ is not attained, then 
$$
\sigma_{\mathrm{pt}}(T_X) =\emptyset;\;	
\sigma_{\mathrm{r}}(T_X)=\mathrm{int}(\mathcal{R}_{q_X}).
$$
whereas if both $q_X\ge2$ and $q_X$ is attained, then 
$$
\sigma_{\mathrm{pt}}(T_X)=\emptyset;\; 
\sigma_{\mathrm{r}}(T_X)=  \mathcal{R}_{q_X}\setminus\{\pm1\}.
$$
\item [(iii)] Suppose that  $q_X\le2\le p_X$ and, for those cases when either $p_X=2$ 
or $q_X=2$ occur,  that they are not attained. Then 
$$
\sigma_{\mathrm{pt}}(T_X) =\sigma_{\mathrm{r}}(T_X)=\emptyset;\; 
\sigma_{\mathrm{c}}(T_X)=\sigma(T_X). 
$$
\end{itemize}
\end{proposition}

Not all cases are covered by the previous result. The identification of 
$\sigma_{\mathrm{c}}(T_X)$ encounters serious difficulties. More precise information,
which we now present, is possible when it is known that 
$\underline{\alpha}_X = \overline{\alpha}_X$. Observe that the union of two sets of the
form $\mathcal{R}_s$, for $1<s<\infty$, is again a set of the same form, clearly 
the larger one. Using this fact and an interpolation argument (via Boyd's theorem)
it can be shown, \cite[Proposition 7.1]{COR-am}, that
\begin{equation}\label{8.7}
\sigma(T_{X})\subseteq 
\mathcal{R}_{1/\underline{\alpha}_X}\cup \mathcal{R}_{1/\overline{\alpha}_X},
\end{equation}
whenever $X$ is a separable r.i.\ space on $(-1,1)$ with 
$0<\underline{\alpha}_X\le \overline{\alpha}_X<1$. The previous containment,
together with earlier results, allows a complete identification of the fine spectra
of $T_X$ whenever $0<\underline{\alpha}_X = \overline{\alpha}_X<1$,
in which case $p_X=q_X=1/\underline{\alpha}_X=1/\overline{\alpha}_X$.
This identification is provided in \cite[Theorem 7.2]{COR-am} which we formulate
in the following table. The issue of whether the indices $p_X, q_X$ are attained or not
attained is indicated by a.\ or n.a., respectively.

\bigskip
\begin{center}
\begin{tabular}{ |c | c| c | c |  c| }
\hline
	$\sigma(T_X)=\mathcal{R}_{p_X}$  & a./n.a. &
	$\sigma_{\mathrm{pt}}(T_X)$& $\sigma_{\mathrm{r}}(T_X)$
	&$\sigma_{\mathrm{c}} (T_X)$
\\ \hline
	$p_X<2$ &n.a. &$\mathrm{int}(\mathcal{R}_{p_X})$ & 
	$\emptyset$ &$\partial\mathcal{R}_{p_X}	$
\\ \hline
   	& a.  & $\mathcal{R}_{p_X}\setminus\{\pm1\}$&$\emptyset$  
	&$\{\pm1\}$
\\ \hline
  	$p_X>2$ & n.a. &$\emptyset$&$\mathrm{int}(\mathcal{R}_{p_X})$ &
  	$\partial\mathcal{R}_{p_X}$
\\ \hline
	&a. &$\emptyset$&$\mathcal{R}_{p_X}\setminus\{\pm1\}$ &
  	$\{\pm1\}$
\\ \hline
	$p_X=2$ & $p_X$ a.  &$(-1,1)$& $\emptyset$& $\{\pm1\}$
\\ \hline
	& $q_X$ a. &$\emptyset$& $(-1,1)$& $\{\pm1\}$
\\ \hline
	& $p_X,q_X$ n.a.  &$\emptyset$& $\emptyset$& $\mathcal{R}_2=[-1,1]$
   	\\ \hline
		
\end{tabular}
\end{center}
\bigskip

Supplementary to the ``standard'' known r.i.\ spaces $X$ for which
$0<\underline{\alpha}_X = \overline{\alpha}_X<1$ is satisfied, some additional examples
of such spaces are  studied in  \cite[\S7]{COR-am}. For instance,
this is the case for the family of all separable Orlicz spaces $L^\Phi(-1,1)$  whose 
Young function $\Phi$ and its complementary Young function both satisfy the 
$\Delta_2$-condition, and $\Phi$ satisfies the condition
$$
\lim_{t\to0^+}\frac{t\Phi'(t)}{\Phi(t)}=\lim_{t\to\infty}\frac{t\Phi'(t)}{\Phi(t)},
$$
\cite[Theorem 1.3]{FK}. The same is true for the classical Lorentz space 
$\Lambda^p(w)$ on $(-1,1)$, for $1\le p<\infty$ with 
$w$ a positive, decreasing and continuous
function on $(0,2)$ satisfying   $\lim_{t\to0^+}w(t)=\infty$. 
Not so well known, perhaps, is the grand Lebesgue space $L^{p)}$, 
for $1<p<\infty$,  introduced in \cite{IS}. Its associate space $(L^{p)})'$, 
a so-called small  
Lebesgue space,    is a separable, non-reflexive r.i.\ space
with both Boyd indices  equal to $1/p$, \cite[Theorem 2.1]{FG}.

It is also possible to find a smaller superset (namely, $\mathcal{R}_{p_X}$)
of $\sigma(T_X)$ than that given in \eqref{8.7} by using the indices $p_X, q_X$ in place of
$\underline{\alpha}_X, \overline{\alpha}_X$ and applying a different interpolation theorem
than that of Boyd, namely one for interpolation spaces between
$L^p$ and $L^q$, for $1<p<q<\infty$. These have a description as spaces which are
both interpolation spaces between $L^1$ and $L^q$ and between $L^p$ and $L^\infty$,
\cite{AC}, \cite{LS}. This approach yields our final result concerning 
$\sigma(T_X)$; see  Theorems 7.7 and 7.9 and Propositions 7.6 and 7.8 in \cite{COR-am}. 
Note that $\sigma_{\mathrm{c}}(T_X)$ is precisely identified.

\begin{theorem}\label{t-8.10}
Let $X$ be a separable r.i.\ space on $(-1,1)$ with non-trivial Boyd indices.
\begin{itemize}
\item[(i)] Suppose that $2>p_X=q_X$ and $q_X$ is attained. If  $X$ is 
an interpolation space between $L^2$ and $L^{p_X}$, then
$$
\sigma_{\mathrm{pt}}(T_X) =\mathrm{int}(\mathcal{R}_{p_X});
\quad
\sigma_{\mathrm{r}}(T_X)=\emptyset;
\quad
\sigma_{\mathrm{c}}(T_X)=\partial\mathcal{R}_{p_X}.
$$
\item[(ii)] Suppose that $p_X=q_X>2$ and $p_X$ is attained. 
If $X$ is an interpolation space 
between $L^2$ and $L^{p_X}$, then
$$
\sigma_{\mathrm{pt}}(T_X) =\emptyset;
\quad
\sigma_{\mathrm{r}}(T_X)=\mathrm{int}(\mathcal{R}_{p_X});
\quad
\sigma_{\mathrm{c}}(T_X)=\partial\mathcal{R}_{p_X}.
$$
\end{itemize}
\end{theorem}


\section{Integral representation of the finite Hilbert transform}
\label{S9}

Let $(\Omega,\Sigma,\nu)$ be a finite  measure space (always positive) and $X(\nu)$
be a B.f.s.\ over this measure space which contains $L^\infty(\nu)$.
Given a Banach space $Y$, a continuous linear operator $T\colon X(\nu)\to Y$
generates a finitely additive $Y$-valued measure $m_T\colon\Sigma\to Y$ via
\begin{equation}\label{9.1}
m_T(A):=T(\chi_A),\quad A\in\Sigma.
\end{equation}
For each $\Sigma$-simple function $s=\sum_{j=1}^na_j\chi_{A_j}$ with
$\{a_j\}_{j=1}^n\subseteq\C$ and $\{A_j\}_{j=1}^n\subseteq\Sigma$
one can define the integrals
$$
\int_A s\,dm_T:=\sum_{j=1}^na_jm_T(A\cap A_j),\quad A\in\Sigma.
$$
Observe, if $A\in\Sigma$ is a $\nu$-null set, then it is also a 
$m_T$-\textit{null set}, meaning that $m_T(B)=0$ for every $B\subseteq A$
with $B\in\Sigma$. In this generality not much can be said about   
the interaction between the properties of $m_T$ and those of $T$. 
However, if $m_T$ is actually
$\sigma$-additive on the $\sigma$-algebra $\Sigma$ (which is assumed
to be the case henceforth and is automatic whenever $X(\nu)$ has
a.c.-norm), then one can define the space $L^1_w(m_T)$ of all
\textit{scalarly $m_T$-integrable functions}. Namely, it consists of 
those $\Sigma$-measurable functions $f\colon\Omega\to\C$ 
such that $\int_\Omega|f|\,d|\langle m_T,y^*\rangle|<\infty$ for each $y^*\in Y^*$, where 
$\langle m_T,y^*\rangle$ denotes the complex measure $A\mapsto 
\langle m_T(A),y^*\rangle$ on $\Sigma$ and $|\langle m_T,y^*\rangle|$
is its  variation measure. The space $L^1_w(m_T)$ is a Banach space for the norm
$$
\|f\|_{L^1_w(m_T)}:=\sup_{\|y^*\|\le1}\int_\Omega|f|\,d|\langle m_T,y^*\rangle|,
\quad f\in L^1_w(m_T).
$$
There is an important closed ideal $L^1(m_T)$ of $L^1_w(m_T)$ which
consists of those functions $f\in L^1_w(m_T)$ with the additional property that,
for each $A\in\Sigma$, there exists an element $\int_Af\,dm_T\in Y$ 
(necessarily unique) satisfying
$$
\Big\langle \int_Af\,dm_T,y^*\Big\rangle
=\int_Af\,d\langle m_T,y^*\rangle,\quad y^*\in Y^*.
$$
Elements of $L^1(m_T)$ are called $m_T$-\textit{integrable functions}.
Clearly the space $\text{sim}(\Sigma)$ of all $\Sigma$-simple functions is
a vector subspace of $L^1(m_T)$. The \textit{integration map 
associated with $m_T$} is the linear operator $I_{m_T}\colon L^1(m_T)\to Y$ defined by
\begin{equation}\label{9.2}
I_{m_T}(f):=\int_\Omega f\,dm_T,\quad f\in L^1(m_T);
\end{equation}
it satisfies $\|I_{m_T}\|=1$. The restriction of the norm $\|\cdot\|_{L^1_w(m_T)}$
to $L^1(m_T)$ is denoted by $\|\cdot\|_{L^1(m_T)}$. A finite measure
$\lambda\colon\Sigma\to[0,\infty)$ is called a \textit{control measure}
for $m_T$ if $\lambda$ and $m_T$ have the same null sets. In this case
$L_w^1(m_T)$ is a B.f.s.\ over $(\Omega,\Sigma,\lambda)$.
Moreover, $L^1(m_T)$ always has a.c.-norm but, it may  fail the Fatou property. 
Rybakov's theorem, 
\cite[Theorem IX.2.2]{DU}, asserts that there exists $y^*\in Y^*$ such that 
$|\langle m_T,y^*\rangle|$ is a control measure for $m_T$. Such a vector 
$y^*$ is called a \textit{Rybakov functional} for $m_T$. The operator
$T$ is called $\nu$-\textit{determined} if $\nu$ is a control measure for $m_T$.
In this case there is an intimate connection between $m_T$ and $T$ as
seen in the following result, \cite[Theorem 4.14]{ORS-P}.

\begin{proposition}\label{p-9.1}
Let $(\Omega,\Sigma,\nu)$ be a finite measure space and $X(\nu)$ 
be a B.f.s.\  over   $(\Omega,\Sigma,\nu)$ with a.c.-norm. 
Let $Y$ be a Banach space and $T\colon X(\nu)\to Y$ be a continuous 
linear operator which is $\nu$-determined. Then $L^1(m_T)$ is the largest amongst all B.f.s.'
(over $(\Omega,\Sigma,\nu)$) with a.c.-norm such that 
$X(\nu)\subseteq L^1(m_T)$, with the natural inclusion being continuous, 
and to which $T$ admits a $Y$-valued continuous linear extension.

Such an extension is unique and is precisely the integration map
$I_{m_T}\colon L^1(m_T)\to Y$ given in \eqref{9.2}, that is
\begin{equation}\label{9.3}
T(f)=\int_\Omega f\,dm_T,\quad f\in X(\nu).
\end{equation}
\end{proposition}

All of the above facts about vector measures, their associated integration 
maps   and their $L^1$-spaces can be found in \cite{C-PhD}, \cite{DU}, \cite{KK},
\cite{Le-1}, \cite{Le-2}, \cite{ORS-P}, for example, and the references therein. 
The space $L^1(m_T)$ is also called the \textit{optimal (a.c.-) domain}
of $T$. This ``optimal extension process'' has been investigated for kernel operators
(cf.\  \cite{AS}, \cite{CR(i)},  \cite{Sz}), Sobolev imbeddings 
(cf.\ \cite{CR(ii)}, \cite{CR(iii)}, \cite{EKP},  \cite{KP}), the Hardy operator,
\cite{DS}, and the Hausdorff-Young inequality, \cite{MR}.
For convolutions with measures in $L^p$-spaces see  \cite{OR(i)},  \cite{OR(ii)},
\cite[Ch.7]{ORS-P} and for more general Fourier $p$-multiplier
operators we refer to \cite{MOR}.

As will be seen, it is the integral representation of $T$ given by \eqref{9.3} 
which allows the use of the well developed theory of integration with respect 
to vector measures to deduce operator theoretic properties of the FHT operator $T$. 
We will now concentrate our attention on the case when $T$ is the FHT on $(-1,1)$
given in \eqref{T}. Whenever $X$  is a r.i.\ space
on $(-1,1)$ with non-trivial Boyd indices
the associated vector measure $m_{T_X}$  specified by \eqref{9.1}
is denoted simply by $m_X$. Let now $(\Omega,\Sigma,\nu)$ be the 
particular measure space $((-1,1), \mathcal{B},\mu)$, where $\mu$ denotes Lebesgue
measure (\textit{only} for this section). 
A crucial fact is that $T_X$ is $\mu$-determined, \cite[Proposition 3.2(iv)]{COR-qm}.

The question arises of whether there exists an optimal domain 
for $T_X$ beyond the class of B.f.s.' \textit{having a.c.-norm}. 
This point was addressed in Proposition \ref{p-6.3}, where it was explained that
$[T_X,X]$ is the largest B.f.s.\  on $(-1,1)$ containing $X$ continuously relative to the norm
$\|\cdot\|_{[T_X,X]}$
and to which $T_X\colon X\to X$ has a continuous $X$-valued extension.
Moreover, $[T_X,X]$ has the Fatou property. It was this space $[T_X,X]$ 
which was used to show that $T_X$, for $X$ r.i.\ 
with non-trivial Boyd indices, is already 
optimally defined; see Theorem \ref{t-6.1}. So, how are the two optimal domains spaces,
$L^1(m_X)$ with a.c.-norm and $[T_X,X]$ with the Fatou property,
related? The clue lies in the special role played by the closed ideal
$X_a\subseteq X$. Lemma \ref{l-4.1} implies that $L^q\subseteq X_a$ for all $q$ satisfying
$L^q\subseteq X$; see the proof of Lemma 2.3 in \cite{COR-qm}. In particular, 
$X_a\not=\{0\}$. Note that $X_a$ may \textit{not} have the Fatou
property (cf.\ \cite[Remark 3.12(b-2)]{COR-qm}). Since the range
$$
m_X(\mathcal{B})=\{T_X(\chi_A):A\in \mathcal{B}\}\subseteq X_a
$$
whenever $X$ has  non-trivial Boyd indices, 
\cite[Lemma 2.3]{COR-qm}, it follows that $m_X$ is always 
$\sigma$-additive when interpreted as an $X_a$-valued vector measure and hence,
also as an $X$-valued vector measure. 
Moreover, $T_X(X_a)\subseteq X_a$, \cite[Lemma 2.3]{COR-qm}.
This allows the possibility for the following refinement of Proposition 
\ref{p-9.1}, for the case when the operator $T$ there is replaced with $T_X$;
see \cite[Lemma 2.6]{COR-qm}.

\begin{proposition}\label{p-9.2}
Let $X$ be a r.i.\ space on $(-1,1)$ with non-trivial Boyd indices. Then $L^1(m_X)$ 
is the largest amongst all B.f.s.'   on $(-1,1)$ with a.c.-norm 
into which $X_a\not=\{0\}$ is continuously embedded and to 
which the restriction $T|_{X_a}$ admits
an $X$-valued continuous linear extension. Further, such an extension is unique
and equals the integration map $I_{m_X}\colon L^1(m_X)\to X$.

In particular, $X_a\subseteq L^1(m_X)$ continuously and 
$I_{m_X}(f)=T|_{X_a}(f)\in X_a\subseteq X$ for all $f\in X_a$.
\end{proposition}

So, we see that necessarily  $L^1(m_X)\subseteq [T_X,X]$ and
$[T_X,X]$ has the Fatou property. It was noted that the B.f.s.\ 
$L^1_w(m_X)$ also has the Fatou property and satisfies $L^1(m_X)\subseteq L^1_w(m_X)$.
It is time to make the precise connection between the various spaces
involved. The following result is a combination of Lemma 3.9, Theorem
3.10 and Corollary 3.11 in \cite{COR-qm}. See also Corollary \ref{c-6.4}.

\begin{proposition}\label{p-9.3}
Let $X$ be a r.i.\ space on $(-1,1)$ with non-trivial Boyd indices. 
The following statements are valid.
\begin{itemize}
\item[(i)] $L^1_w(m_X)$ has the Fatou property and coincides with the
bi-associate space $L^1(m_X)''$ of $L^1(m_X)$.
\item[(ii)] The following conditions are equivalent.
\begin{itemize}
\item[(a)] $L^1(m_X)$ has the Fatou property.
\item[(b)] $L_w^1(m_X)$ has a.c.-norm.
\item[(c)] $L^1(m_X)=L_w^1(m_X)$.
\item[(d)] $\mathrm{sim}(\mathcal{B})$ is dense in $L_w^1(m_X)$.
\end{itemize}
\item[(iii)] The B.f.s.\ $L_w^1(m_X)$ is the minimal B.f.s.\ on 
$(-1,1)$ with the Fatou property which contains (with norm $\le1$) $L^1(m_X)$.
\item[(iv)] The natural inclusions
$$
X_a\subseteq L^1(m_X)\subseteq [T_X,X]=X\subseteq L^1
$$
hold and are continuous. Furthermore,
$$
\|f\|_{L^1(m_X)}=\|f\|_{[T_X,X]},\quad f\in L^1(m_X).
$$
In particular, $L^1(m_X)$ is a closed ideal in $[T_X,X]=X$.

\noindent
Moreover, we have the integral representation of $T_X$ given by
\begin{equation}\label{9.4}
T_X(f)=I_{m_X}(f)=\int_{(-1,1)}f\,dm_X,\quad f\in L^1(m_X).
\end{equation}
\item[(v)] $L^1(m_X)=X_a$ with equivalent lattice norms.
\item[(vi)] $L_w^1(m_X)=[T_X,X]$ identically as B.f.s.' on $(-1,1)$.
\item[(vii)] If $X$ has a.c.-norm, then $X= L^1(m_X)=[T_X,X]$.
\item[(viii)] The natural inclusion $X\subseteq L^1(m_X)$ holds if and only if
$T_X(X)\subseteq X_a$.
\end{itemize}
\end{proposition}

We now record various properties of  $m_X$ itself. The \textit{variation measure} $|m_X|$  of the
vector measure $m_X$ is defined as for scalar measures,
\cite[Definition I.1.4]{DU}, by replacing the absolute value with
the norm $\|\cdot\|_X$. A subset $C$ of $X$ is called 
\textit{order bounded} if there exists $0\le f\in X$ such that $|h|\le f$ for all 
$h\in C$. The following result is Proposition 3.2 of  \cite{COR-qm}.
We point out that there is a typing error in the identity (3.4) in the
proof of part (iii) of Proposition 3.2 in \cite{COR-qm}; the correct formula is
$$
g_0:=\widecheck{T}_{L^{p'}}(\sigma)=-wT(\sigma/w).
$$

\begin{proposition}\label{p-9.4}
Let $X$ be a r.i.\ space on $(-1,1)$ with non-trivial Boyd indices. 
\begin{itemize}
\item[(i)] $m_X(\mathcal{B})\subseteq X_a$.
\item[(ii)] For every $g\in X'\subseteq X^*$ the complex measure 
$\langle m_X,g\rangle$ is given by
$$
\langle m_X,g\rangle(A)=-\int_AT_{X'}(g)\,d\mu,\quad A\in\mathcal{B}.
$$
Since $T_{X'}(g)\in X'\subseteq L^1$, the variation measure 
of $\langle m_X,g\rangle$ is given by
$$
|\langle m_X,g\rangle|(A)=\int_A|T_{X'}(g)|\,d\mu,\quad A\in\mathcal{B}.
$$
\item[(iii)] There exists a Rybakov functional $g_0\in X'$ 
satisfying $|\langle m_X,g_0\rangle|=\mu$. In particular, the $\mu$-null
and the $m_X$-null sets coincide.
\item[(iv)] The variation $|m_X|(A)=\infty$ for every non-$m_X$-null set
$A\in\mathcal{B}$. 
\item[(v)] $m_X(\mathcal{B})$ is not a relatively compact subset of $X$.
\item[(vi)] $m_X(\mathcal{B})$ is not an order bounded subset of $X$.
\end{itemize}
\end{proposition}

It should be noted that the range of every $\sigma$-additive vector measure defined on a 
$\sigma$-algebra is necessarily relatively weakly compact, \cite[Theorem 2.9]{BDS}.

The integral representation \eqref{9.4} for $T_X$, together with the special features of 
$m_X$ listed in Proposition \ref{p-9.4}, provide
a means to deduce various operator theoretic properties of $T_X$. 
This is illustrated by Corollary 3.4 of  \cite{COR-qm}, where it is established 
that the operator $T_X\colon X\to X$ is neither order bounded
(as it maps the order bounded subset $\{\chi_A:A\in\mathcal{B}\}$
of $X$ to the non-order bounded subset $m_X(\mathcal{B})$
of $X$), nor is it completely continuous, nor is it compact. The arguments used 
rely on the principle that certain properties of a vector 
measure are closely related  to the membership of its integration map in appropriate
operator ideals; see, for example, \cite{ORR-P(i)},
\cite{ORR-P(ii)}, \cite{ORR-P(iii)} and the references therein.
As a sample, if $I_{T_X}$ was compact, then $m_X$ necessarily 
has finite variation,  \cite[Theorem 4]{ORR-P(i)}. Since this is not the 
case (cf.\ Proposition \ref{p-9.4}(iv)) and $T_X=I_{m_X}$ whenever 
$X$ has a.c.-norm (which implies that
$X=L^1(m_X)$ via 
Proposition  \ref{p-9.3}(vii)),  it follows that $T_X$ is not a compact operator.
Alternatively, since $[-1,1]$ is an uncountable set, this also follows from
Proposition \ref{p-8.3}(c).

We end this section with a brief discussion of the case when $X=L^1$. 
Here the indices $\underline{\alpha}_{L^1}=\overline{\alpha}_{L^1}=1$ are trivial
and so the above results are not applicable as $T(L^1)\not\subseteq L^1$. Nevertheless, since $m_{L^2}\colon\mathcal{B}\to L^2$ is 
$\sigma$-additive and the natural inclusion $j\colon L^2\to L^1$ 
is continuous, the set function 
$m_{L^1}:=j\circ m_{L^2}\colon\mathcal{B}\to L^1$ is a 
$\sigma$-additive vector measure.
Moreover, since $L^1$ is weakly sequentially complete (hence, it 
cannot contain a copy of $c_0$), it is known that
$L^1(m_{L^1})=L^1_w(m_{L^1})$, \cite[Theorem 5.1]{Le-2}. Noting that
$T(s)=\int_{(-1,1)}s\,dm_{L^1}\in L^1$, for each
$s\in\mathrm{sim}(\mathcal{B})$, and that $\mathrm{sim}(\mathcal{B})$ 
is dense in $L^1(m_{L^1})$, \cite[Theorem 3.7(ii)]{ORS-P},
it follows that the integration map $I_{m_{L^1}}\colon
L^1(m_T)\to L^1$ satisfies
\begin{equation}\label{9.5}
I_{m_{L^1}}(f)=\int_{(-1,1)}f\,dm_{L^1}=T(f),\quad f\in L^1(m_{L^1})\subseteq L^1.
\end{equation}
We have seen in Section \ref{S7} that
$T(\logl)\not\subseteq \logl$ but, $T$ does map $\logl$
continuously into the strictly larger space $L^1$ (this operator
was denoted by $\tlog$) with no further extension possible.
It turns out, somewhat
remarkably, that $L^1(m_{L^1})=\logl$ with equivalent lattice norms 
so that \eqref{9.5} becomes
$$
\tlog(f)=\int_{(-1,1)}f\,dm_{L^1},\quad f\in \logl.
$$
Moreover, $m_{L^1}$ and $\mu$ have the same null sets and there 
exists a Rybakov functional $g_0\in L^\infty=(L^1)^*=(L^1)'$ for
$m_{L^1}$ satisfying $\mu=|\langle m_{L^1},g_0\rangle|$. So,
$\tlog$ is surely $\mu$-determined. The proof of the existence of
$g_0$ for $m_{L^1}$ is significantly more involved than for $m_{X}$
and is based on an analysis of $T$ acting on certain H\"older continuous functions.
Since $L^\infty\subseteq X'\subseteq X^*$ for all r.i.\ spaces $X$ on $(-1,1)$ 
and $m_X(A)=m_{L^1}(A)$ for all $A\in\mathcal{B}$, this is a considerable
strengthening of Proposition \ref{p-9.4}(iii)  because the \textit{same function}
$g_0\in L^\infty$ can be chosen as a Rybakov functional satisfying
$\mu=|\langle m_{X},g_0\rangle|$ for \textit{every} such $X$.
Moreover, the vector measure $m_{L^1}\colon\mathcal{B}\to L^1$ has infinite 
variation over every set $A\in\mathcal{B}$ satisfying $\mu(A)>0$ and,
for such a set $A$, the subset $\{m_{L^1}(A\cap B): B\in\mathcal{B}\}$
is not order bounded in $L^1$. It is also the case that 
$m_{L^1}(\mathcal{B})$ is not a relatively compact subset of $L^1$.
For all of the above facts we refer to Section 3 of \cite{COR-prep}.
Using these properties of $m_{L^1}$, together with $I_{m_{L^1}}=\tlog$, it 
is established in Section 4 of  \cite{COR-prep} that the operator $\tlog$
is not order bounded, not completely continuous (hence, not compact) 
and also not weakly compact.



\end{document}